\documentclass[journal,twoside,web]{ieeecolor}

\let\labelindent\relax
\usepackage{generic}
\usepackage{cite}
\usepackage{amsmath,amssymb,amsfonts}
\usepackage{algorithmic}
\usepackage{graphicx}
\usepackage{textcomp}
\usepackage{comment}

\usepackage{hyperref}%
\hypersetup{colorlinks=true, linkcolor=blue, breaklinks=true, urlcolor=blue}
\usepackage{doi}

\usepackage[all]{xy}
\usepackage{amsmath}
\usepackage{amscd}
\usepackage{amsthm}
\usepackage{mathrsfs}
\usepackage{amssymb}
\usepackage{tikz}
\usepackage[utf8]{inputenc}
\usepackage[yyyymmdd]{datetime}
\usepackage{graphicx}
\usepackage{subcaption}
\usepackage{lipsum}
\usepackage{multicol}
\usepackage{mathtools}
\usepackage{multirow}
\usepackage{threeparttable}
\usepackage{pdfpages}
\usepackage{multicol} 
\usepackage{mathtools}
\usepackage{xspace}
\usepackage{url,doi}
\usepackage{amsmath}

\newcommand{\subscr}[2]{#1_{\textup{#2}}}
 \newcommand{\setdef}[2]{\{#1
	\; | \; #2\}}

\newcommand{\map}[3]{#1: #2 \rightarrow #3}

\newcommand{\suchthat}{\;\ifnum\currentgrouptype=16 \middle\fi|\;}
\newcommand{\scirc}{\raise1pt\hbox{$\,\scriptstyle\circ\,$}}

\newcommand{\real}{\mathbb{R}}
\newcommand{\realnonnegative}{\mathbb{R}_{\geq0}}

\DeclareMathOperator{\diag}{diag}
\DeclareMathOperator{\spn}{span}

\DeclareMathOperator{\sgn}{sgn}
\renewcommand{\top}{\mathsf{T}}  

\DeclareSymbolFont{bbold}{U}{bbold}{m}{n}
\DeclareSymbolFontAlphabet{\mathbbold}{bbold}
\newcommand{\vect}[1]{\mathbbold{#1}}

\newcommand\oprocendsymbol{\hbox{$\triangle$}}
\newcommand\oprocend{\relax\ifmmode\else\unskip\hfill\fi\oprocendsymbol}

\usepackage[prependcaption,colorinlistoftodos]{todonotes}
\newcommand{\todoinb}[1]{\todo[inline,color=gnblue6!10, linecolor=gnblue6!10]{\small#1}}
\newcommand{\todoing}[1]{\todo[inline,color=green!20, linecolor=orange!250]{\small#1}}

\newcommand{\fbtodo}[1]{\todoinb{\textbf{FBTodo from FB:} #1}}
\newcommand{\kevintodo}[1]{\todoinb{\textbf{Todo from Kevin:} #1}}

\usepackage{enumitem} \renewcommand{\theenumi}{(\roman{enumi})}

\newcommand{\mcK}{\mathcal{K}}

\usepackage{xcolor}
\usepackage{xparse}
\usepackage{mathrsfs}
\usepackage{ stmaryrd }
\usepackage{amssymb}
\usepackage{graphics} 
\usepackage{amsmath} 
\usepackage{amssymb}  
\usepackage{amsthm}
\usepackage{enumitem} 

\usepackage[prependcaption, colorinlistoftodos]{todonotes}

\newcommand{\gdptodo}[1]{\todoing{\textbf{From GDP:} #1}}

\DeclareMathOperator*{\sign}{sign}

\renewcommand{\natural}{{\mathbb{N}}}

\newcommand{\integernonnegative}{\ensuremath{\mathbb{Z}}_{\ge 0}}

\newcommand{\eps}{\epsilon}
\newcommand{\rhoess}{\subscr{\rho}{ess}}

\renewcommand{\theenumi}{(\roman{enumi})}

\DeclareSymbolFont{bbold}{U}{bbold}{m}{n}
\DeclareSymbolFontAlphabet{\mathbbold}{bbold}
\usepackage{soul}

\NewDocumentCommand\seminorm{mg}{\left\vert\kern-0.25ex\left\vert\kern-0.25ex\left\vert #1 \right\vert\kern-0.25ex\right\vert\kern-0.25ex\right\vert^{\mcK} \IfNoValueF{#2}{_{#2}}}

\NewDocumentCommand\seminormnok{mg}{\left\vert\kern-0.25ex\left\vert\kern-0.25ex\left\vert #1 \right\vert\kern-0.25ex\right\vert\kern-0.25ex\right\vert \IfNoValueF{#2}{_{#2}}}
\usepackage{centernot}
\usepackage{wrapfig}
\usepackage{cancel}[Smaller]

\NewDocumentCommand\bigseminorm{mg}{\Big\lvert\!\Big\lvert\!\Big\lvert #1 \Big\rvert\!\Big\rvert\!\Big\rvert^\mcK \IfNoValueF{#2}{_{#2}}}

\DeclareMathOperator{\conv}{conv}

\newtheorem{theorem}{Theorem}
\newtheorem{remark}[theorem]{Remark}
\newtheorem{example}[theorem]{Example}
\newtheorem{lemma}[theorem]{Lemma}
\newtheorem{definition}[theorem]{Definition}

\newtheorem{corollary}[theorem]{Corollary}
\newtheorem{conjecture}[theorem]{Conjecture}

\NewDocumentCommand\norm{mg}{\lVert #1 \rVert \IfNoValueF{#2}{_{#2}}}
\definecolor{gnred}{RGB}{255,91,89}

\definecolor{gnred1}{RGB}{71,0,0} 
\definecolor{gnred2}{RGB}{117,0,0} 
\definecolor{gnred3}{RGB}{164,0,0} 
\definecolor{gnred4}{RGB}{211,0,0} 
\definecolor{gnred5}{RGB}{255,0,0} 
\definecolor{gnred6}{RGB}{255,42,34} 
\definecolor{gnred7}{RGB}{255,91,89} 

\definecolor{gnblue1}{RGB}{0,36,71}   
\definecolor{gnblue2}{RGB}{0,60,118}  
\definecolor{gnblue3}{RGB}{0,85,164}  
\definecolor{gnblue4}{RGB}{0,108,212} 
\definecolor{gnblue5}{RGB}{0,133,255}  
\definecolor{gnblue6}{RGB}{35,156,255} 
\definecolor{gnblue7}{RGB}{88,177,255} 

\definecolor{gnbrown1}{RGB}{71,27,0}  
\definecolor{gnbrown2}{RGB}{117,45,0} 
\definecolor{gnbrown3}{RGB}{164,62,0} 
\definecolor{gnbrown4}{RGB}{211,80,0} 
\definecolor{gnbrown5}{RGB}{255,97,0} 
\definecolor{gnbrown6}{RGB}{255,127,26} 
\definecolor{gnbrown7}{RGB}{255,155,86} 

\newcommand{\dist}{\operatorname{dist}}
\newcommand{\proj}{\operatorname{proj}}

\usepackage{version}

\includeversion{arxiv}
\excludeversion{tac}

\usepackage{layouts}
\def\BibTeX{{\rm B\kern-.05em{\sc i\kern-.025em b}\kern-.08em
    T\kern-.1667em\lower.7ex\hbox{E}\kern-.125emX}}
\begin{document}

\title{Dual Seminorms, Ergodic Coefficients \\ and Semicontraction Theory}
\author{Giulia De Pasquale \IEEEmembership{Student Member, IEEE}, Kevin D. Smith \IEEEmembership{Student Member, IEEE},  Francesco Bullo \IEEEmembership{Fellow, IEEE}, M. Elena Valcher \IEEEmembership{Fellow, IEEE}
\thanks{Submitted on \today. This work was supported in part by AFOSR grant FA9550-22-1-0059 and by Fondazione Ing.\ Aldo Gini. }
\thanks{Giulia De Pasquale and Maria Elena Valcher are with the Dipartimento di Ingegneria dell'Informazione Università di Padova, Padova, 35131,  Italy.{\tt\small \{giulia.depasquale, meme\}@dei.unipd.it.
}}
\thanks{Kevin D. Smith and Francesco Bullo are with the Center for Control, Dynamical Systems, and Computation, UC Santa Barbara, Santa Barbara, CA 93101 USA. {\tt\small \{kevinsmith, bullo\}@ucsb.edu.
}}}
\maketitle
\begin{abstract}
  Dynamical systems that are contracting on a subspace are said to be
  semicontracting.  Semicontraction theory is a useful tool in the study of
  consensus algorithms and dynamical flow systems such as Markov chains.

  To develop a comprehensive theory of semicontracting systems, we
  investigate seminorms on vector spaces and define two canonical notions:
  projection and distance seminorms. We show that the well-known $\ell_p$
  ergodic coefficients are induced matrix seminorms and play a central role
  in stability problems. In particular, we formulate a duality theorem that
  explains why the Markov-Dobrushin coefficient is the rate of contraction
  for both averaging and conservation flows in discrete time. Moreover, we
  obtain parallel results for induced matrix log seminorms.  Finally, we
  propose comprehensive theorems for strong semicontractivity of linear and
  non-linear time-varying dynamical systems with invariance and
  conservation properties both in discrete and continuous time.
\end{abstract}

\section{Introduction}

\subsection*{Problem description and motivation}
Before Stefan Banach proved his famous contraction principle in 1922~\cite{SB:1922}, Andrey Markov started in 1906~\cite{AAM:1906} the study of stochastic processes. As documented by Eugene Seneta~\cite{ES:06-Markov}, Markov established a key contraction inequality and a corresponding contraction factor now known with the name of \emph{ergodic coefficient} of a Markov chain. This paper aims to provide a modern semicontraction theory approach to explain and generalize ergodic coefficients.

To be concrete, let the matrix $A$ be row-stochastic and consider the discrete-time dynamical systems
\begin{subequations}\label{eq:A+Atopdt}
 \begin{alignat}{2}
 x(k+1) &= Ax(k),  \label{eq:A+Atopdt-A}
 \\
 \pi(k+1) &= A^\top\pi(k). \label{eq:A+Atopdt-AT}
 \end{alignat}
\end{subequations}
Similarly, let $L$ be a Laplacian matrix and consider the continuous-time
counterparts:
\begin{equation}\label{eq:A+Atopct}
  \dot{x}(t)=-L x(t), \quad \dot{\pi}(t)= -L^\top \pi(t).
\end{equation}
These systems are perhaps the simplest examples of general averaging-based dynamics (e.g., robotic coordination and distributed optimization) and dynamical flow systems (e.g., compartmental and traffic systems). Important generalizations include systems of the form $\dot{x}=f(t,x)$, where $f$ satisfies invariance properties (generalizing $A\vect{1}_n=\vect{1}_n$) or conservation properties (generalizing $\vect{1}_n^\top A^\top=\vect{1}_n^\top$); in all these (linear and nonlinear) cases, the system 
is at most marginally stable.

Markov and later scientists essentially showed that, under a certain connectivity assumption, maps of the form $\pi \mapsto A^\top\pi$ are contraction maps with respect to the total variation distance on the simplex. To be specific, define the simplex $\Delta_n=\setdef{x\in\real^n}{x\geq0, \vect{1}_n^\top{x}=1}$ and the total variation distance on $\Delta_n$ by $d_{{\rm TV}}(\pi,\sigma)= \tfrac{1}{2} \sum_i|\pi_i-\sigma_i|$. Then any two solutions $\pi(k),\sigma(k)$ to~\eqref{eq:A+Atopdt-AT} satisfy
\begin{equation} \label{eq:Markov-bound}
    d_{{\rm TV}}\big(\pi(k)-\sigma(k)\big) \leq \tau_1(A)^k d_{{\rm TV}}\big(\pi(0)-\sigma(0)\big),
\end{equation}
where $\tau_1(A)$ is the so-called \emph{Markov-Dobrushin ergodic coefficient} defined by
\begin{equation}
\label{def:tau1}
    \tau_1(A) := \max_{\norm{z}{1} = 1, \, z^\top \vect{1}_n =0} \norm{A^\top z}{1}.
\end{equation}
In short, when $\tau_1(A)<1$, existence, uniqueness and global exponential stability of an equilibrium $\pi^*\in\Delta_n$  for system~\eqref{eq:A+Atopdt-AT} is ensured.

Now comes a remarkable similarity. If one defines the seminorm  $\seminormnok{x}{\dist,\infty} = \tfrac{1}{2} \left( \max_i\{x_i\} - \min_j\{x_j\} \right)$, the following fact is also known~\cite[Theorem~1.1]{DJH:98} about averaging systems of the form~\eqref{eq:A+Atopdt-A}:
\begin{equation} \label{eq:averaging-bound}
    \seminormnok{x(k)}{\dist,\infty} \leq \tau_1(A)^k  \seminormnok{x(0)}{\dist,\infty}.
\end{equation}

Despite the extensive research in this field, numerous known related facts remain somehow mysterious and numerous related mathematical questions remain open. 
For example, why is the same ergodic coefficient $\tau_1$ relevant for the contraction properties of both dynamical flow systems and averaging systems? And is it the tightest such bound?
How does one generalize the bounds~\eqref{eq:Markov-bound} and~\eqref{eq:averaging-bound} to ergodic coefficients $\tau_p$ defined with respect to arbitrary $\ell_p$ norms (instead of the $\ell_1$ norm in~\eqref{def:tau1})?
How does one provide a unified robust stability analysis for both systems?
What are the canonical Lyapunov functions for both systems \eqref{eq:A+Atopdt-A}-\eqref{eq:A+Atopdt-AT}, whose discrete-time variation along the flow is described by $\tau_p(A)$?  
How does one define ergodic coefficients for continuous-time systems?
Is there a contraction theoretic framework that applies to time-varying and nonlinear systems with generalized invariance or conservation properties?

\subsection*{Contributions}
This paper provides a comprehensive answer to all the open research
questions outlined above.

In order to define Lyapunov functions for averaging, flow systems and their
generalizations to nonlinear dynamical systems with invariant subspaces, we
study seminorms, induced matrix seminorms for discrete time systems and
logarithmic seminorms for continuous-time systems.  A key contribution of
this paper is to explain precisely in what sense ergodic coefficients are
\emph{induced matrix seminorms} and, when less than unity,  contraction
factors for discrete-time systems. This equality is the fundamental reason
why ergodic coefficients play a critical role in robust stability theory
for discrete-time dynamical systems with invariance properties.  It is
surprising that induced norms are widely studied in the matrix theory
literature, but induced seminorms much less (e.g., see~\cite{RAH-CRJ:85}).

After characterizing various seminorms' properties, we define two canonical
sets of seminorms, namely, \emph{distance} and \emph{projection seminorms},
and establish remarkable duality properties between the two. Our first
result generalizes and strengthens the so called \emph{Markov contraction
  inequality} as a duality result between the aforementioned seminorms.
Our duality result precisely explains why the induced matrix seminorms for
both $A$ and $A^\top$ are identical, when computed with respect to dual
seminorms.  Particular emphasis is given to the case of \emph{consensus
  seminorms}, that is, seminorms whose kernel is the consensus space (i.e.,
seminorms that are positive definite about the consensus space). Consensus
seminorms appear naturally in averaging algorithms and surprisingly in
systems with conservation property (such as Markov chains and dynamical
flow systems).

It is an elementary algebraic observation that the total variation distance
on the simplex arises from the restriction of the $\ell_1$ projection
consensus seminorm.

We then leverage all these notions to provide a general nonlinear semicontraction theory, grounded in two key theorems both for continuous and discrete time varying dynamical systems.  The semicontraction theory we develop is tailored to systems with invariance or
conservation properties.  More in detail, when either the system's Jacobian leaves
invariant the seminorm kernel (invariance property) or its orthogonal complement  (conservation property), there is a well
defined notion of \emph{perpendicular dynamics} which is strictly
contracting. For both systems, in the linear time varying case, we  show how canonical Lyapunov functions (some of which partly known in the literature) naturally arise from seminorms. For the non-linear case, our first key theorem establishes  conditions and features of strong semicontracting  continuous time, time varying systems that enjoy the invariance property. The theorem extends Theorem 13 in \cite{SJ-PCV-FB:19q} through the formulation of a cascade decomposition and by establishing a strong contractivity property on the orthogonal complement to the seminorm kernel. The second key theorem is entirely novel and  pertains semincontraction conditions  for continuous time, time varying, dynamical systems that enjoy the conservation property. A discrete time version of these two theorems is also provided.


\subsection*{Literature review}

Interest in contractivity of dynamical systems via matrix measures can be traced back to Demidovi\v{c} \cite{BPD:61} and Krasovski\u{\i} \cite{NNK:63}. Logarithmic norms have been exploited in control theory later on by Desoer and Vidyasagar in \cite{CAD-MV:1975} and applied in the study of contraction theory for dynamical systems for the first time by Lohmiller and Slotine \cite{WL-JJES:98}.  
In the context of control theory, this literature inspired many  generalizations of contraction theory such as partial contraction \cite{WW-JJES:05}, weak- and semi-contraction \cite{SJ-PCV-FB:19q}, horizontal contraction on Riemannian and Finsler manifolds \cite{JWSP-FB:12za,FF-RS:14},  etc.

In particular, partial contraction refers to convergence of systems trajectories to a specific behavior,  or a manifold \cite{JJES:03}, see also \cite{MdB-DF-GR-FS:16} for a survey on this theory. While partial contraction establishes convergence to a manifold, semicontraction ensures contractivity on the subspace perpendicular to the kernel of the seminorm. For a characterization of partial contraction in the $\ell_2$-norm for the study of synchronization in networked systems, see \cite{WW-JJES:05}. The notion of partial contraction is closely related to the one of semicontraction and weak contraction proposed and investigated in \cite{SJ-PCV-FB:19q}. Semicontraction theory relies on a relaxed concept of matrix measure, known  as matrix semimeasure. For this reason, contractivity of a dynamical system is only ensured on a certain subspace and the distance between trajectories is allowed to increase along certain directions.

A relevant behavior, to which semicontraction theory applies, is the one of
consensus for dynamical systems. Strictly related to consensus, when it
comes to stochastic systems, is the concept of (weak) ergodicity
\cite{ATS-AJ:08}. The concept of weak ergodicity was first formalized in
1931 by Kolmogorov \cite{AnK:1931}, who stated that a sequence of
stochastic matrices is \emph{weakly ergodic} if the rows of the matrix
product {\color{black} tend to become identical} as the number of factors
increases.  The study of ergodicity coefficients is traced back to the
pioneering work of Markov \cite{AAM:1906}, in 1906, in which a first
expression of ergodicity coefficients was provided in the context of the
Weak Law of Large Numbers. Subsequent works from Doeblin \cite{WD:1937} and
Dobrushin \cite{RLD:1956} provided conditions for weak ergodicity.  The key
results in this research area were extended and then reviewed by Seneta in
the 80's, see, e.g.,~\cite{ES:81}.  A survey of ergodicity coefficients is
given by Ipsen and Selee~\cite{ICFI-TMS:11}. a historical discussion is
given by Hartfiel~\cite[Chapter~1]{DJH:98}, and a recent treatment on their
connection with spectral graph theory is given by Marsli and Hall
\cite{RM-FJH:20}.  A characterization of ``convergability"
\cite{JL-SM-ASM-BDOA-CY:11}, namely the convergence of a product of an
infinite number of stochastic matrices, is studied by Liu et. al in
\cite{JL-SM-ASM-BDOA-CY:11}, where a different approach, based on optimally
deflated matrices, is proposed. Despite the evident relation between
ergodicity coefficients, contraction factors and induced matrix seminorms,
especially in the context of stochastic and averaging systems
\cite{ZA-RF-AH-YC-TTG:20}, to the best of our knowledge none in the past
has shed full light on their connections (see \cite{CB:13a-arxiv} for some
preliminary work in this direction). This manuscript aims to bridge the
existing gap in the scientific literature between semicontraction and
ergodicity of dynamical systems.

\subsection*{Paper organization}
Section \ref{sec:notation} presents notation and preliminary
results. Section \ref{sec:duality} introduces the projection and distance
seminorms and establishes their duality relationship.  Section
\ref{sec:induced_matr} pertains with induced matrix seminorms and induced
matrix log-seminorms. In Section \ref{sec:semicontraction} semicontraction
theory is applied to dynamical systems.  Finally, Section
\ref{sec:conclusions} concludes the manuscript.

All theorems in this manuscript are new. Lemmas and Corollaries are either
new or simple derivations from known results. \begin{tac}Due to page
  constraints some proofs are omitted and reported in the extended
  technical report \cite{GDP-KDS-FB-MEV:21m-arxiv}.\end{tac} \begin{arxiv}
  This manuscript extends the submitted version in the IEEE Transaction on
  Automatic Control and includes the proofs of Lemma~\ref{the:wellp}
  and Theorem~\ref{the:duallog}, explicit expressions for projection seminorms of
  columns stochastic matrices
  \eqref{eq:ind-norm-pi1}-\eqref{eq:ind-norm-piinf}, Corollary
  \ref{thm:mu-explicit}, the property~\ref{propmatr3} from Lemma
  \ref{lem:propertiesIMS}, and Remark~\ref{mylem:normseminorm}. \end{arxiv}
 


\section{Notation and Preliminaries}\label{sec:notation}
\subsection{Notation}
The set $\realnonnegative$ is the set of nonnegative real numbers. Let $I_n\in \real^{n\times n}$ denote the identity matrix of size $n$. Let $\vect{1}_n$ and
$\vect{0}_n$ denote the  $n$ dimensional column vectors all whose entries equal $1$ and $0$, 
respectively.  Let $\vect{e}_i$ denote the $i$-th vector of the canonical basis in $\real^n$. 
For a matrix $A\in \real^{n\times n}$, 
let 
{\color{black}$A^\top$ denote its transpose, $[A]_{i, j}$
 its $(i, j)$th entry. 
 }
The matrix $A$ is \emph{nonnegative} if all its entries
are nonnegative, it is \emph{row stochastic} if it is nonnegative and  $A\vect{1}_n=\vect{1}_n$, it is \emph{column stochastic} if $A^\top$ is row stochastic. 

Given $A\in \real^{n\times{n}}$,  a vector subspace ${\mathcal K} \subseteq \real^n$ is  $A$-{\em invariant} if $A{\mathcal K}\subseteq {\mathcal K}$. The symbol $\map{\langle \cdot, \cdot \rangle }{\real^n \times \real^n}{\real}$ denotes the standard inner product on $\real^n$.
We let $\Pi_\perp$ denote the orthogonal projection matrix onto $\mcK^\perp$,  where the symbol $\mcK^\perp$ denotes the orthogonal complement of $\mcK$. Note that {\color{black} $\Pi_\perp= \Pi_\perp^\top$,} and if $\mcK = {\rm span} \{\vect{1}_n$\}, then 
${\color{black} \Pi_\perp =I_n-\vect{1}_n\vect{1}_n^\top/n =: \Pi_n.}$
 Given $x\in \real^n$, the perpendicular and parallel components of $x$ to $\mcK$ are denoted by $x_\perp = \Pi_\perp x$ and  $x_\parallel = (I_n-\Pi_\perp)x$, respectively. 
Define the \emph{$n$-simplex} {\color{black} as} $\Delta_n = \setdef{v \in \real_{\ge 0}^n}{\vect{1}_n^\top v =1}$ and the \emph{sign function}, $\map{\sign}{\real}{\{-1, 0, 1\}}$, {\color{black} as} $\sign(x) = \frac{x}{|x|}$
if $x \neq 0$,  and $\sign(0)=0$.
 Given two matrices $A,B \in \real^{n\times n}$ we use the notation $A\preceq B$ to indicate that $A-B$ is a negative semidefinite matrix. 

 A directed, weighted graph is a triple \cite{FB:22}, $\mathcal{G} = (\mathcal{V}, \mathcal{E}, \mathcal{A})$, where $\mathcal{V} = \{1,\dots, n\}$ is the set of vertices, $\mathcal{E} \subseteq \mathcal{V} \times \mathcal{V}$ is the set of arcs and $\mathcal{A}$ is the adjacency matrix. An arc $(i,j)$ belongs to $\mathcal{G}$ if and only if $[\mathcal{A}]_{ij} \neq 0$. Two nodes $i,j \in \mathcal{V}$ are \emph{weakly adjacent} if either $(i,j) \in \mathcal{E}$ or $(j,i) \in \mathcal{E}$.

Given a real vector space $V$, the \emph{dual space} $V^\star$ is the vector space of linear maps from $V$ into $\real$. If $V = \real^n$, then $V^\star$ is the vector space of row vectors in $\real^n$. In this case, it is typical to make a slight abuse of notation and assume $V^\star = \real^n$.

\subsection{Basic concepts}
We start with some basic useful concepts.  For $x\in\real^n$ and
$p\in\natural$, the \emph{$\ell_p$-norm} of $x$ is
\begin{equation*}
    \norm{x}{p}\triangleq \Big( \sum_{i=1}^n \lvert x_i\rvert^p
    \Big)^{\frac{1}{p}},
\end{equation*}
while the $\ell_\infty$-norm is
\begin{equation*}
  \norm{x}{\infty} =
  \lim_{p\to \infty} \Big( \sum_{i=1}^n \lvert x_i\rvert^p \Big)^{\frac{1}{p}} = \max_{i} \lvert x_i \rvert. 
\end{equation*}
For $A\in \real^{n \times m}$ and $p \in \natural$, the $\ell_p$-\emph{induced norm} of $A$ is
\begin{equation*}
\norm{A}{p} = \max_{\substack{x\in \real^{m}\\\norm{x}{p}\leq 1}} \norm{Ax}{p}.
\end{equation*}

\begin{definition}[Seminorms]
  A function $\map{\seminormnok{ \cdot} }{\real^n}{\realnonnegative}$  is a \emph{seminorm} on $\real^n$ if it satisfies the following properties for all $x, y \in \real^n$ and $a \in \real$:
\begin{align*}
   & \text{(homogeneity): } \seminormnok{ a x }= \lvert a
  \rvert\seminormnok{ x },\; \text{and}\\ 
  & \text{(subadditivity): } \seminormnok{ x + y}
  \leq \seminormnok{ x}  + \seminormnok{ y}.
\end{align*}
Consequently, $\seminormnok{x} = 0$ does not imply $x = 0$. The \emph{kernel} of a seminorm is the vector space
\[
    \mcK \triangleq \ker(\seminormnok{\cdot})=\left\{x \in \real^n : \seminormnok{x} = 0\right\}.
\] 
\end{definition}
From now onward, for a seminorm $\seminormnok{\cdot}$ on $\real^n$ with kernel $\mcK$ we will use the symbol $\seminorm{\cdot}$.
\begin{lemma}[Seminorms of orthogonal projections] \label{lem:projection}
    Let $\seminorm{\cdot}{}$ be a seminorm on $\real^n$ with kernel $\mathcal{K}$, and let $\Pi_\perp$ be the orthogonal projection matrix onto $\mcK^\perp$. For all $x\in \real^n$, $\seminorm{x} = \seminorm{\Pi_\perp x}$.
\end{lemma}
\begin{proof}
The result is a direct consequence of the reverse triangle inequality and the sub additivity property of seminorms applied to the orthogonal decomposition $x = x_\perp  + x_\parallel $, with $x_\parallel  \in \mathcal{K}$.
\end{proof}
\begin{arxiv}
\begin{remark}[Relationship between norm and seminorm]\label{mylem:normseminorm}
A seminorm on $\real^n$ with kernel $\mcK$ induces a norm on $\mcK^\perp$ by restriction. Vice-versa, given a subspace $\mcK$ of $\real^n$, a norm $\norm{\cdot}$ on  $\mcK^\perp$, denoted by $\norm{\cdot}{\perp}$, induces a seminorm on $\real^n$ with kernel $\mcK$ by projection: $\seminorm{x} = \norm{\Pi_\perp x}{\perp}$. 
\end{remark}
\end{arxiv}


\begin{definition}\emph{(Induced seminorm \cite{SJ-PCV-FB:19q})}
 \label{def:induced_sem}
Given   a seminorm  $\map{\seminorm{\cdot}}{\real^n}{\realnonnegative}$ with kernel $\mathcal{K}$, the induced seminorm on $\real^{n \times n}$ is
\begin{equation*}\label{eq:weight_semi}
    \seminorm{ A} \triangleq \max_{\substack{\seminorm{ x} \leq 1 \\ x \perp\mathcal{K}}} \seminorm{A x}.
\end{equation*}
\end{definition}

\begin{definition}[Matrix logarithmic seminorms \cite{FB:22-CTDS}] 
 Given a seminorm $\map{\seminorm{\cdot}}{\real^n}{\realnonnegative}$ with kernel $\mathcal{K}$, the induced matrix logarithmic seminorm on $\real^{n \times n}$ is
\begin{equation*}
    \mu_{\seminorm{\cdot}}(A) \triangleq \lim_{h\rightarrow0^+} \frac{\seminorm{I_n+hA}-1}{h}.
    \label{eq:log-norm-limit}
\end{equation*}
\end{definition}




\begin{definition}[Generalized $\ell_p$ ergodicity coefficient \cite{UGR-CPT:85}]\label{mydef:erg_coeff}
  Given $p\in [1,\infty]$ and a vector subspace $\mcK \subset\real^m$, the
  \emph{generalized $\ell_p$ ergodicity coefficient}
  $\map{\tau_p}{\mcK\times\real^{m\times{n}}}{\real}_{\geq 0}$ is defined by
  \begin{equation}\label{eq:erg_coeff}
    \tau_p(\mcK,A) := \max_{\substack{\norm{z}{p} = 1 \\ z \perp \mcK}} \norm{A^\top z}{p}.
  \end{equation}
\end{definition}

The ergodicity coefficient \eqref{eq:erg_coeff} is the norm of the operator
defined on the real (normed) linear space $\mcK^\perp$ by $x \rightarrow
A^\top x$ \cite{UGR:75}.




\begin{lemma}[$\ell_2$-Norm LMI characterization \cite{FB:22-CTDS}] \label{lem:lmi}
    Given any $A \in \real^{n \times n}$, 
    \begin{equation*} \label{eq:lmi-2norm}
        ||A||_2 = \min \setdef{b \in \realnonnegative}{
            A^\top A \preceq b^2 I_n} .
    \end{equation*} 
\end{lemma}

\section{Seminorms and Duality}\label{sec:duality}

\subsection{Projection and Distance Seminorms}

\begin{figure*}
    \centering
%
    \includegraphics[width=0.25\linewidth]{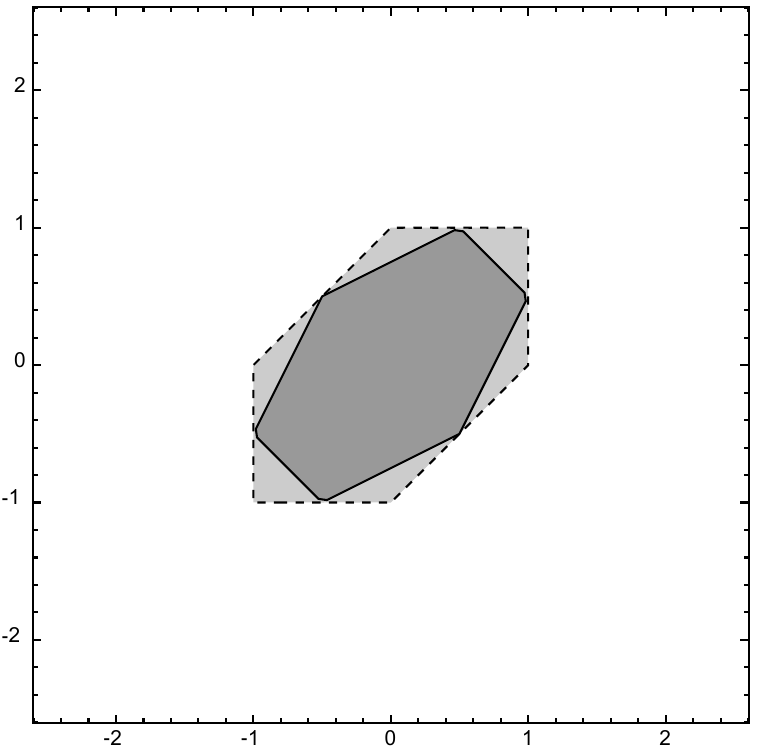} \hfill
    \includegraphics[width=0.25\linewidth]{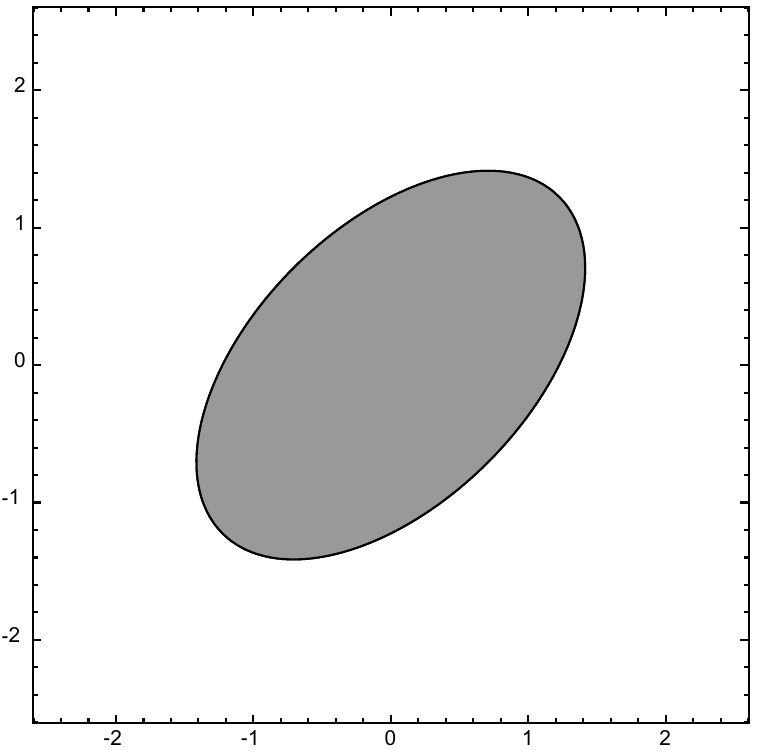} \hfill
    \includegraphics[width=0.25\linewidth]{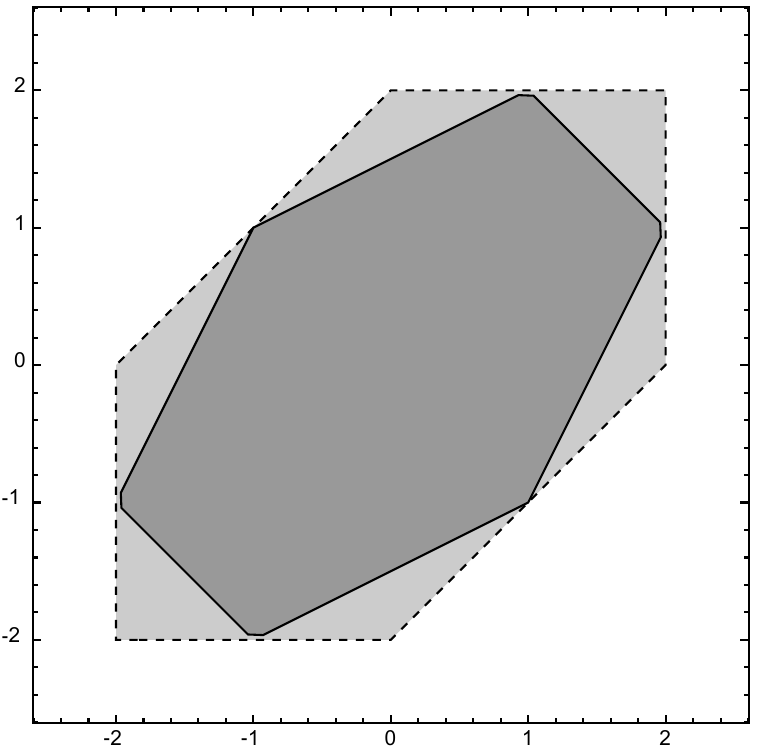} \hfill
    \caption{Two-dimensional sections of three-dimensional unit disks of projection (solid contours) and distance (dashed contours) consensus seminorms. We plot the section corresponding to $(x_1, x_2, x_3=0)$ for $p=1$ (left), $p=2$ (center), and $p=\infty$ (right).}
    \label{fig:sublevel-sets}
\end{figure*}
In the following we provide the definition of projection and distance seminorms. These two  seminorms will play a fundamental role in the duality result.
\begin{definition}[Projection and distance seminorms]
    Let $\mcK \subset \real^n$ be a vector space and $\Pi_\perp \in
    \real^{n \times n}$ be the orthogonal projection matrix onto
    $\mcK^\perp$. For each $p\in[1,\infty]$, define the
    \emph{$\ell_p$-projection seminorm} with respect to $\mcK$ by
    \begin{equation} \label{eq:def-proj}
      \seminorm{x}{\proj, p} \triangleq \norm{\Pi_\perp x}{p} 
    \end{equation}
    and the \emph{$\ell_p$-distance seminorm} with respect to $\mcK$ by
    \begin{equation} \label{eq:def-dist}
        \seminorm{x}{\dist, p} \triangleq \dist_p(x, \mcK) = \min_{u \in \mcK} \norm{x-u}{p}.
    \end{equation}
\end{definition}
Note that the optimization problem \eqref{eq:def-dist} is well posed since
the norm function is convex.
\begin{lemma}[Basic properties] \label{lem:seminorm-props}
  For each $p \in [1, \infty]$,
  \begin{enumerate}
  \item \label{prop:kernel} $\ker(\seminorm{\cdot}{\proj, p}) = \ker(\seminorm{\cdot}{\dist, p}) = \mcK$,
  \item \label{prop:order-1} $\seminorm{x}{\dist, p} \le \min\{\norm{x}{p}, \seminorm{x}{\proj, p}\}$ for all $x \in \real^n$.
  \end{enumerate}
\end{lemma} 
\begin{proof}
  Statement \ref{prop:kernel} is obvious from \eqref{eq:def-proj} and
  \eqref{eq:def-dist}. Next, we compute
  \begin{align*}
    \min_{u \in \mcK} \norm{x - u}{p} &\le \norm{x - \vect 0_n}{p} = \norm{x}{p},\\
    \min_{u \in \mcK} \norm{x - u}{p} &\le \norm{x - (I_n - \Pi_\perp)x}{p} = \seminorm{x}{\proj, p}.
  \end{align*}
  This completes the proof of statement~\ref{prop:order-1}.
\end{proof}
It is not true in general that $\seminorm{x}{\proj, p} \le \norm{x}{p}$.


\begin{example}[Seminorms for consensus and stationary distribution]
    When $\mcK = \spn\{\vect 1_n\}$ and $p \in \{1, 2, \infty\}$, explicit
    formulas for the $\ell_p$-projection and distance seminorms are either
    easily derivable or available in the literature
    \cite{MH-LVDH:84,JL-SM-ASM-BDOA-CY:11,FB:22-CTDS}. For each $x \in
    \real^n$, with the shorthand $x_{\rm avg} = \frac 1 n \vect 1_n^\top
    x$,
    \begin{align*}
        \seminorm{x}{\proj,1} &= \sum_{i=1}^n \left|x_i - x_{\rm avg}\right|, \\
        \seminorm{x}{\proj,2} &=  \Big( \frac{1}{n} \sum\nolimits_{i, j} (x_i - x_j)^2 \Big)^{1/2}, \\
        \seminorm{x}{\proj,\infty} &= \max_{i} \left|x_i - x_{\rm avg} \right|.
    \end{align*}
    Next, sort the entries of $x$ according to $x_{(1)} \ge x_{(2)} \ge
    \cdots \ge x_{(n)}$. With this shorthand,
    \begin{align*}
      \seminorm{x}{\dist,1} &= \sum_{i = 1}^{\lfloor \frac n 2 \rfloor} x_{(i)} - \sum_{i = \lceil \frac n 2 \rceil + 1}^n x_{(i)}, \\
      \seminorm{x}{\dist,2} &=  \Big( \frac{1}{n} \sum\nolimits_{i, j} (x_i - x_j)^2 \Big)^{1/2}, \\
      \seminorm{x}{\dist,\infty} &= \frac{1}{2} \left( x_{(1)} - x_{(n)} \right)
      =
      \frac{1}{2} \big( \max_i\{x_i\} - \min_j\{x_j\} \big).
    \end{align*}   
    Figure \ref{fig:sublevel-sets} illustrates the unit disks for these
    seminorms on~$\real^3$.
\end{example}

\begin{example}[Total variation distance]
    The \emph{total variation} \cite[Section~4.1]{DaL-YP:17} is a metric on the simplex $\Delta_n$ defined by
    \[
        d_{{\rm TV}}(x,y) \triangleq \frac{1}{2}\sum_{i=1}^n |x_i-y_i|.
    \] 
    Given any two vectors $x, y \in \Delta_n$, a simple derivation shows
    \[
        d_{{\rm TV}} (x,y)= \frac{1}{2}\seminorm{x-y}{\proj,1},
    \] 
    where $\seminorm{\cdot}{\proj, 1}$ is the $\ell_1$-projection seminorm with respect to the kernel $\mcK = \spn\{\vect 1_n\}$.
\end{example}

\subsection{Duality}

In this section we establish a useful duality relationship between
projection and distance seminorms. We start with the notion of dual
seminorm.

%

\begin{definition}[Dual seminorm]\label{def:dseminorm}
  Let $\seminorm{\cdot}$ be a seminorm on a real vector space $V\subseteq
  \real^n$ with kernel $\mcK \subset V$. The \emph{dual seminorm} is the
  function $\map{\seminorm{\cdot}_\star}{V^\star}{\real}$ defined by
  \begin{equation*}
    \seminorm{x}{\star} \triangleq \max_{\substack{\seminorm{y} \leq 1\\ y \perp \mathcal{K}}} \langle x, y \rangle.
  \end{equation*}
\end{definition}
\begin{tac}We omit the proof of the following natural result.\end{tac}
\begin{lemma}[Well-posedness of dual seminorms]\label{the:wellp}
  Let $\seminorm{\cdot}{}$ be a seminorm on a real vector space $V$ with
  kernel $\mcK$. Then the dual seminorm $\seminorm{\cdot}_\star$ is a
  seminorm on $V^\star$.
\end{lemma}
\begin{arxiv}
\begin{proof}
    Let $y, z \in V^\star$, and let $a \in \real$. Since $x = \vect 0_n$ satisfies $\seminorm{x} \le 1$ and $x \in \mcK^\perp$, $\seminorm{y}{\star} \ge \langle y, \vect 0_n \rangle = 0$, establishing the non-negativity of $\seminorm{\cdot}{\star}$. To prove homogeneity,
    \begin{align*}
        \seminorm{a y}{\star} &= \max_{\substack{\seminorm{x} \le 1 \\ x \in \mcK^\perp}} \langle ay, x \rangle
        = \max_{\substack{\seminorm{x} \le 1 \\ x \in \mcK^\perp}} |a| \sgn(a) \langle y, x \rangle \\
        &= |a| \max_{\substack{\seminorm{x} \le 1 \\ x \in \mcK^\perp}} \langle y, \sgn(a) x \rangle = |a| \seminorm{y}{\star}.
    \end{align*}
    Finally, to prove sub-additivity,
    \begin{align*}
        \seminorm{y + z}{\star} &= \max_{\substack{\seminorm{x} \le 1 \\ x \in \mcK^\perp}} \langle y + z, x \rangle
        = \max_{\substack{\seminorm{x} \le 1 \\ x \in \mcK^\perp}} \langle y, x \rangle + \langle z, x \rangle \\
        &\le \max_{\substack{\seminorm{x} \le 1 \\ x \in \mcK^\perp}} \langle y, x \rangle + \max_{\substack{\seminorm{x} \le 1 \\ x \in \mcK^\perp}} \langle z, x \rangle 
        = \seminorm{y}{\star} + \seminorm{z}{\star}.
    \end{align*}
\end{proof} 
\end{arxiv}

When $V=\real^n$, we make the usual identification $(\real^n)^\star =
\real^n$. In this case, the kernel of the dual seminorm is identical to the
kernel of the primal seminorm.

Next, we present an important generalization to arbitrary $\ell_p/\ell_q$
norms of the Markov contraction inequality from \cite[Lemma
  2.3]{MH-LVDH:84}.
\begin{lemma}[Markov contraction inequality]\label{mci}
  Let $p,q\in[1,\infty]$ satisfy $p^{-1} + q^{-1} = 1$ (with the convention
  $1/\infty=0$) and consider a vector space $\mcK \subset \real^n$. For all
  $x, y \in \real^n$,
  \begin{equation*} \label{lem:markov}
    x^\top \Pi_\perp y \le \seminorm{x}{\proj,p} \seminorm{y}{\dist,q}. 
  \end{equation*}
\end{lemma}

\begin{proof}
  For each $u \in \mcK$ satisfying $\seminorm{y}{\dist, q} = \norm{y -
    u}{q}$, 
  \begin{align*}
    x^\top \Pi_\perp y = x^\top \Pi_\perp(y - u) 
    & \overset{\text{(H\"older's ineq)}}{\le} \norm{\Pi_\perp x}{p} \norm{y - u}{q}.
  \end{align*}
  The result follows from minimizing with respect to $u$.
\end{proof}

\begin{remark}[Markov contraction and H\"older's inequalities]
  For the inner product of vectors perpendicular to a subspace, the Markov
  contraction inequality provides a tighter bound than the H\"older's
  inequality $x^\top y \leq \norm{x}{p} \norm{y}{q}$. In fact, as a
  consequence of Lemma~\ref{lem:seminorm-props}\ref{prop:order-1},
  \begin{equation}
    x^\top \Pi_\perp y \le \seminorm{x}{\proj,p} \seminorm{y}{\dist,q} \leq \seminorm{x}{\proj,p}\seminorm{y}{\proj,q}. \notag
  \end{equation}
\end{remark} 

Next, we recall that, for unconstrained vectors, the H\"older's inequality
provides a tight bound in the sense that, for all $x \in \real^n$, there
exists $y \in \real^n$ such that $x^\top y = \norm{x}{p} \norm{y}{q}$.  We
now show this tightness result also for the Markov contraction inequality,
thereby establishing the duality relationship between projection and the
distance seminorms.

\begin{theorem}[Duality of distance and projection seminorms]\label{dualitytheorem}
  Let $p,q \in [1,\infty]$ satisfy that $p^{-1}+q^{-1} = 1$ (with the
  convention $1/\infty=0$) and let $\mcK \subset \real^n$ be a vector
  subspace. Then $\seminorm{\cdot}{\dist,p}$ and
  $\seminorm{\cdot}{\proj,q}$, with kernel $\mcK$, are dual seminorms:
    \begin{align} 
        \seminorm{\cdot}{\dist, p} &= \left( \seminorm{\cdot}{\proj, q} \right)_\star \label{eq:firstdual_g}\\
        \seminorm{\cdot}{\proj, q} &= \left( \seminorm{\cdot}{\dist, p} \right)_\star. \label{eq:seconddual_g}
    \end{align}
\end{theorem}

\begin{proof}
  To prove \eqref{eq:firstdual_g}, consider two cases. First, if $x \in
  \mcK$, then $\seminorm{x}{\dist,p} = 0$.  On the other hand, $y \in
  \mcK^\perp$ implies $y^\top x = 0$. So, both sides of
  \eqref{eq:firstdual_g} are zero.  Second, if $x \not\in \mcK$, by Lemma
  \ref{lem:coefs} in Appendix~\ref{sec:appendix}, there exists $\psi_p(x)
  \in \mcK^\perp$ with $\seminorm{\psi_p(x)}{\proj,q} = 1$ such that
  \[
  \seminorm{x}{\dist,p} = \psi_p(x)^\top x \le \max_{\substack{\seminorm{y}{\proj,q} \le 1 \\ y \in \mcK^\perp}} y^\top x = \left( \seminorm{\cdot}{\proj, q} \right)_\star.
  \]
  To prove the opposite inequality, choose any $y \in \real^n$ such that
  $\seminorm{y}{\proj,q} \le 1$ and $y \in \mcK^\perp $. Then $||y||_q =
  \seminorm{y}{\proj,q} \le 1$, so by Lemma \ref{mci},
  \[
  y^\top x \le ||y||_q \seminorm{x}{\dist,p} \le \seminorm{x}{\dist,p}
  \] 
  
  To prove equality \eqref{eq:seconddual_g} we notice, as in the previous
  first case, that if $x \in\mcK$, then $\seminorm{x}{\proj,q} = 0$, while
  $ y \in \mcK^\perp$ implies that $y^\top x = 0$, so both sides of
  \eqref{eq:seconddual_g} are zero.  Otherwise, in the second case, if $x
  \not\in\mcK$, by Lemma \ref{lem:coefs_disa}, there exists $\zeta_q(x) \in
  \mcK^\perp$ with $\seminorm{\zeta_q(x)}_{\dist,p}\leq 1$ such that
  \[
  \seminorm{x}{\proj,q} = \zeta_q(x)^\top x \le \max_{\substack{\seminorm{y}_{\dist,p} \leq 1 \\ y \in \mcK^\perp}} y^\top x = \left( \seminorm{\cdot}{\dist, p} \right)_\star.
  \]
  To prove the opposite inequality, choose any $y \in \real^n$ such that
  $\seminorm{y}{\dist,p} \le 1$ and $y \in \mcK^\perp$.  Lemma~\ref{mci}
  implies
  \[
  x^\top y = x^\top\Pi_\perp y \le \seminorm{x}{\proj,q} \seminorm{y}{\dist,q} \le  \seminorm{x}{\proj,q}.
  \]
  This concludes the proof.
\end{proof}

\section{Induced Matrix Seminorms and Log Seminorms}\label{sec:induced_matr}


\subsection{Induced Matrix Seminorms}

In the following we list some basic properties related to induced matrix seminorms.
\begin{lemma}[Properties of induced matrix seminorms]\label{lem:propertiesIMS}
Let $\seminorm{\cdot}$ be a seminorm on $\real^n$ with kernel $\mcK$. For any $A,B\in \real^{n \times n}$,
\begin{enumerate}
    \item \label{prop1} $\seminorm{Ax} \leq \seminorm{A} \seminorm{x}$ for all $x \in \mathcal{K}^\perp$.
\end{enumerate}
\begin{arxiv}
If $A^\top \mathcal{K} \subseteq \mathcal{K}$, then
\begin{enumerate}
    \setcounter{enumi}{1}
    \item \label{propmatr3} $\seminorm{BA} \leq \seminorm{B} \seminorm{A}$.
\end{enumerate}
Moreover,
\end{arxiv}
If $A \mathcal{K} \subseteq \mathcal{K}$, then
\begin{enumerate}
\begin{arxiv}
    \setcounter{enumi}{2}
    \end{arxiv}
    \begin{tac}
     \setcounter{enumi}{1}
    \end{tac}
    \item \label{propmatr1} $\seminorm{A} = \max_{\seminorm{x}\leq 1} \seminorm{Ax}$,
    \item \label{propmatr2} $\seminorm{Ax} \leq \seminorm{A} \seminorm{x}$, and
    \item \label{propmatr4}$\seminorm{AB} \leq \seminorm{A}\seminorm{B}$.
\end{enumerate}

\end{lemma}
\begin{proof}
Property \ref{prop1}  was proven  in \cite{VVK:83}.
\begin{arxiv}
To prove property \ref{propmatr3}, assume $x^* \in \mathcal{K}^\perp$  to be a vector  with $\seminorm{x^*} = 1$ such that $\seminorm{BA x^*} = \seminorm{BA}$. Note that, under the assumption $A^\top \mathcal{K} \subseteq \mathcal{K}$, $Ax^* \in \mathcal{K}^\perp$. Hence
\begin{multline}
    \seminorm{BA} = \seminorm{BA x^*} =\\ \bigseminorm{B \frac{Ax^*}{\seminorm{Ax^*}}} \seminorm{Ax^*} \leq \seminorm{B} \seminorm{A}. \notag 
\end{multline}
\end{arxiv}
To prove property \ref{propmatr1},  decompose any vector $x \in \real^n$ as $x = x_\perp  + x_\parallel $, with $x_\perp  \in \mathcal{K}^\perp$ and $ x_\parallel  \in \mathcal{K}$, and  notice that
\begin{align}
    &\max_{\seminorm{x}\leq 1} \seminorm{Ax} = \max_{\seminorm{x}\leq 1} \seminorm{A(x_\perp  + x_\parallel )} \notag \\
    &=\max_{\seminorm{ x_\perp}\leq 1} \seminorm{A x_\perp } = \max_{\substack{\seminorm{y} \leq 1\\ y \perp \mathcal{K}}} \seminorm{Ay} = \seminorm{A}, \notag
\end{align}
where the second equality is based on Lemma \ref{lem:projection} and exploits the fact that $A\mathcal{K} \subseteq \mathcal{K}$. 
To prove property \ref{propmatr2} we notice that, by adopting the same decomposition as before
\begin{align}
    \seminorm{Ax} &= \seminorm{A x_\perp } = \seminorm{A(x_\perp  / \seminorm{x_\perp })} \seminorm{x_\perp } \notag  \\
    &\leq \seminorm{A} \seminorm{x_\perp } = \seminorm{A} \seminorm{x} \notag
\end{align}
where the first equality is based on Lemma \ref{lem:projection} and exploits  the fact that $A\mathcal{K} \subseteq \mathcal{K}$,  while the inequality derives from Definition \ref{def:induced_sem}. Property \ref{propmatr4} can be  found in \cite{VVK:83}\begin{arxiv} and can be proved by following arguments similar to property \ref{propmatr3}\end{arxiv}.
\end{proof}

Based on Theorem \ref{dualitytheorem} we are now in the position to provide
one of the main results of this manuscript.

For a matrix $A\in \real^{k\times n}$, and for $p,q \in
[1,\infty]$, with $p^{-1}+q^{-1} = 1$ it holds that
 \begin{equation}\label{dualmatrix}
 \norm{A}{p} = \norm{A^\top}{q},
 \end{equation}
 where $\map{\norm{\cdot}{p}}{\real^n}{\realnonnegative}$ and $\map{\norm{\cdot}{q}}{\real^n}{\realnonnegative}$ are dual norms. 

The following theorem represents a generalization of the duality relationship between induced matrix norms \eqref{dualmatrix} to seminorms.

\begin{theorem}[Duality of induced matrix seminorms]\label{mainth}
    Let $p,q \in [1,\infty]$ such that $p^{-1}+q^{-1} = 1$. For any matrix $A \in \real^{n \times n}$, and any vector space $\mcK \subseteq \real^n$,
    \begin{equation} \label{eq:ind-equiv-general}
        \seminorm{A^\top}{\proj,q} = \seminorm{A}{\dist,p}.
    \end{equation}
    Additionally, if 
    $A\mcK \subseteq \mcK$, then
    \begin{equation} \label{eq:ind-equiv-stoch}
        \tau_q(\mcK,A) = \seminorm{A^\top}{\proj,q} = \seminorm{A}{\dist,p}.   \end{equation}
\end{theorem}
\begin{proof}
    Eqn. \eqref{eq:ind-equiv-general} is a direct consequence of Theorem \ref{dualitytheorem}: 
    \begin{align*}
        \seminorm{A^\top}{\proj,q}&=\max_{\substack{\seminorm{x}{\proj,q}\leq 1 \\ x \perp \mcK}} \seminorm{A^\top x}{\proj,q}\notag \\ &\overset{\eqref{eq:seconddual_g}}{=}\max_{\substack{\seminorm{x}{\proj,q} \leq 1 \\ x\perp \mcK}}\  \max_{\substack{\seminorm{y}{\dist,p}\leq 1\\ y \perp \mcK}} y^\top A^\top x \notag\\
        &=\max_{\substack{\seminorm{y}{\dist,p}\leq 1\\ y \perp \mcK}}\ \max_{\substack{\seminorm{x}{\proj,q} \leq 1 \\ x\perp \mcK}} x^\top A y\\
        &\overset{\eqref{eq:firstdual_g}}{=} \max_{\substack{\seminorm{y}{\dist,p} \leq 1 \\ y \perp \mcK}} \seminorm{Ay}{\dist,p} = \seminorm{A}{\dist,p}.\notag
    \end{align*}
To prove \eqref{eq:ind-equiv-stoch} note that

\begin{multline}
    \seminorm{A^\top}{\proj,q} = \max_{\substack{ \seminorm{x}{\proj,q} \leq 1 \\ x\perp \mcK}} \seminorm{A^\top x}{\proj,q} \\=  
    \max_{\substack{ \norm{\Pi_\perp x}{q} \leq 1 \\ x \perp \mcK}} \norm{\Pi_\perp A^\top x}{q}= \max_{\substack{ \norm{x}{q} \leq 1 \\ x \perp \mcK}} \norm{ A^\top x}{q} =  \tau_q(\mcK,A) \notag
\end{multline}
where the second-last equality follows from the fact that $A^\top \mcK^\perp \subseteq \mcK^\perp$ and, since $x \in \mcK^\perp$, $x = \Pi_\perp x$. 
\end{proof} 

In the following we provide some explicit expressions for the distance seminorm of row-stochastic matrices \begin{arxiv}and the projection seminorm of column stochastic matrices\end{arxiv} for the case in which the kernel of the seminorms is the consensus subspace. The explicit expressions can be derived by the ones available in the literature for ergodicity coefficients \cite{ES:84,ICFI-TMS:11}\begin{tac}.\end{tac} \begin{arxiv}and by the duality result from Theorem \ref{mainth}. \end{arxiv} \begin{tac} The expression for the projection seminorm of column stochastic matrices can be easily derived by the duality result from Theorem \ref{mainth} and hence omitted.\end{tac}
\begin{corollary}[Formulas for induced matrix seminorms]
Consider the consensus distance
 \begin{arxiv} and  projection \end{arxiv}seminorm.
    Let $A \in \real^{n \times n}$.
    \begin{arxiv}Let $a_{i, (j)}$ represent the entries of each row $i \in \{1, 2, \dots, n\}$ sorted according to $a_{i, (1)} \ge a_{i, (2)} \ge \cdots \ge a_{i, (n)}$  and similarly, let $a_{(i), j}$ represent the entries of each column $j \in \{1, 2, \dots, n\}$ sorted by
        $a_{(1), j} \ge a_{(2), j} \ge \cdots \ge a_{(n)j}.$\end{arxiv}
        \begin{tac}
         Assume that the entries of each column $j \in \{1, 2, \dots, n\}$ are sorted so that $a_{ (1),j} \ge a_{ (2),j} \ge \cdots \ge a_{(n)},j$.
         \end{tac}

    If $A$ is row-stochastic, then
 	\begin{align}
	    \label{eq:ind-norm-d1} \seminorm{A}{\dist,1} &= \max_j \Bigg\{
	        \sum_{i = 1}^{\lfloor \frac n 2 \rfloor} a_{(i), j} - \sum_{i = \lceil \frac n 2 \rceil + 1}^n a_{(i), j}
	    \Bigg\}, \\
	    \label{eq:ind-norm-d2} \seminorm{A}{\dist,2} &= ||\Pi_n A||_2 =  \min \left\{ b \ge 0 : A^\top \Pi_n A \preceq b^2 \Pi_n\right\}, \\
 	    \label{eq:ind-norm-dinf} \seminorm{A}{\dist,\infty} &= \frac 1 2 \max_{i \ne j} \sum_{k = 1}^n |a_{ik} - a_{jk}| \notag\\
 	    &= 1 - \min_{i \ne j} \sum_{k=1}^n \min\{a_{ik}, a_{jk}\}.
	\end{align} 
	\begin{arxiv}
	If $A$ is column-stochastic, then
 	\begin{align}
	    \label{eq:ind-norm-pi1} \seminorm{A}{\proj,1} &= \frac 1 2 \max_{i \ne j} \sum_{k = 1}^n |a_{ki} - a_{kj}| \notag \\
	    &= 1 - \min_{i \ne j} \sum_{k=1}^n \min\{a_{ki}, a_{kj}\}, \\ 
	    \label{eq:ind-norm-pi2} \seminorm{A}{\proj,2} &= ||\Pi_n A||_2\notag \\
	    &= \min 
     \left\{ b \ge 0 : A^\top \Pi_n A \preceq b^2 \Pi_n\right\}, \\
 	    \label{eq:ind-norm-piinf} \seminorm{A}{\proj,\infty} &= \max_i \left\{
	        \sum_{j = 1}^{\lfloor \frac n 2 \rfloor} a_{i, (j)} - \sum_{j = \lceil \frac n 2 \rceil + 1}^n a_{i, (j)}
	    \right\}.
	\end{align} 
	\end{arxiv}
\end{corollary}

\begin{proof}
    The formulas \eqref{eq:ind-norm-d1},\eqref{eq:ind-norm-dinf} and the first equality in \eqref{eq:ind-norm-d2} follow from the  equivalence $\seminorm{A}{\dist,p} = \tau_q(\vect{1}_n,A)$ in Theorem \ref{mainth} and by applying the explicit expressions for $\tau_q(\vect{1}_n,A)$  provided in Theorem 3.7, Corollary 3.9, Theorem 4.2, and Theorem 6.19 from \cite{ICFI-TMS:11}. 
    
    The second equality in \eqref{eq:ind-norm-d2} follows from Lemma \ref{lem:lmi}, since
\begin{align*}
        ||\Pi_n A||_2 &= \min_{b \in \real} \left\{ (\Pi_n A)^\top(\Pi_n A) \preceq b^2 I_n \right\} \notag  \\
        &= \min_{b \in \real} \left\{ A^\top \Pi_n A \preceq b^2 I_n\right\}.
\end{align*}
     Since $b^2 \Pi_n \preceq b^2 I_n$, it is clear that $A^\top \Pi_n A \preceq b^2 \Pi_n$ implies $A^\top \Pi_n A \preceq b^2 I_n$.  Conversely, assume $A^\top \Pi_n A \preceq b^2 I_n$, so that $v^\top A^\top \Pi_n A v \le b^2 v^\top v$ for all $v \in \real^n$. Then for any $u \in \real^n$, we can decompose $u = u_\bot + u_\parallel$, with $u_\bot \in \spn\{\vect 1_n\}^\bot$ and $u_\parallel \in \spn\{\vect 1_n\}$. Since $A$ is row stochastic,
    \[
        u^\top A^\top \Pi_n A u = u_\bot^\top A^\top \Pi_n A u_\bot \le b^2 u_\bot^\top u_\bot = b^2 u^\top \Pi_n u
    \]
    and thus $A^\top \Pi_n A \preceq b^2 \Pi_n$.  This way we have proved that $A^\top \Pi_n A \preceq b^2 I_n$ if and only if $A^\top \Pi_n A \preceq b^2 \Pi_n$. In turn, this implies
    \[
        ||\Pi_n A||_2= \min_{b \in \real} \left\{ A^\top \Pi_n A \preceq b^2 \Pi_n\right\}.
    \]
\begin{arxiv}  Formulas $\eqref{eq:ind-norm-pi1}-\eqref{eq:ind-norm-piinf}$ are derived by duality. \end{arxiv}
\end{proof}

Finally, we include a comparative analysis for induced seminorms and the notion of optimal deflation given by~\cite{JL-SM-ASM-BDOA-CY:11}. 

\begin{definition}[$p$-optimal deflation \cite{JL-SM-ASM-BDOA-CY:11}]\label{def:optdef}
    For each $p \in [1, \infty]$, the \emph{$p$-optimal deflation of a  matrix $A \in \real^{k \times n}$} is 
    \begin{equation}\label{eq:optdef}
        |A|_p \triangleq \min_{v \in \real^n} ||A - \vect 1_n v^\top ||_p.
    \end{equation}
\end{definition}


\begin{lemma}[Bounds on matrix seminorms] \label{prop:inequalities}
  Given a row-stochastic matrix $A \in \real^{n \times n}$, for each $p \in [1, \infty]$
    \begin{equation*}
        \seminorm{A}{\dist,p} \le |A|_p \le ||A||_p.
    \end{equation*}
\end{lemma}

\begin{proof}
  We first establish that $|A|_p \ge \seminorm{A}{\dist,p}$. By the max-min
  inequality \cite[Section~5.4.1]{SB-LV:04},
  \begin{align*}
    |A|_p &= \min_{v \in \real^n} \max_{||w||_p \leq 1} \left|\left|\left( A - \vect 1_n v^\top \right) w\right|\right|_p \\
    &\ge \max_{||w||_p \leq 1} \min_{v \in \real^n} \left|\left| (Aw) - (v^\top w) \vect 1_n \right|\right|_p \\
    &\ge \max_{||w||_p \leq 1} \seminorm{Aw}{\dist,p}.
  \end{align*}
  Let $w \in \real^n$ be such that $\seminorm{w}{\dist,p} \leq 1$, which
  implies that $||w - \alpha \vect 1_n||_p \leq 1$ for some $\alpha \in
  \real$. Let $u = w - \alpha \vect 1_n$, and observe that $||u||_p \leq
  1$, and that $\seminorm{Au}{\dist,p} = \seminorm{Aw}{\dist,p}$, since $A$
  is row-stochastic and $\seminorm{\cdot}{\dist,p}$ is invariant with
  respect to perturbations in $\spn\{\vect 1_n\}$. Therefore
  \begin{multline*}
    \max_{||w||_p \leq 1} \seminorm{Aw}{\dist,p} \ge \max_{\seminorm{w}{\dist,p} \leq 1} \seminorm{Aw}{\dist,p}  \\
    \ge \max_{\substack{\seminorm{w}{\dist,p} \leq 1 \\ \vect 1_n^\top w = 0}} \seminorm{Aw}{\dist,p} = \seminorm{A}{\dist,p}.
  \end{multline*}
  The inequality $|A|_p \le ||A||_p$ is obtained at $v=\vect 0_n$
  in~\eqref{eq:optdef}.
\end{proof}

\subsection{Induced Matrix Log Seminorms}
 
We now present a duality result for induced matrix log seminorms which is parallel to the one in Theorem \ref{mainth}.
\begin{theorem}[Dual logarithmic seminorms]\label{the:duallog} Let $p,q \in [1,\infty]$  be such that $p^{-1}+q^{-1} =1$. For any matrix $M \in \real^{n \times n}$, and any kernel $\mcK$,
\begin{equation*}
    \mu_{\dist,p} (M) = \mu_{\proj,q}(M^\top).
\end{equation*}
\end{theorem}
\begin{tac}
\begin{proof}
The equality directly follows from the duality of distance and projection induced matrix seminorms.
\end{proof}
\end{tac}
\begin{arxiv}
\begin{proof}
    From \eqref{eq:log-norm-limit} and Theorem \ref{mainth},
\begin{align*}
        &\mu_{\dist,p}(M) = \lim_{h \to 0^+} \frac{\seminorm{I_n + h M}{\dist,p} - 1}{h}=\\
        &\lim_{h \to 0^+} \frac{\seminorm{I_n + h M^\top}{\proj,q} - 1}{h} 
        = \mu_{\proj,q}(M^\top). \notag 
\end{align*}
\end{proof}
\end{arxiv}
We derive now explicit formulas for $\ell_p$-distance logarithmic seminorm
of (minus) Laplacian matrices, for $p \in \{1,2,\infty\}$.
\begin{theorem}[Explicit formulas for distance logarithmic seminorms] \label{thm:mu-explicit}
  Consider the consensus distance and projection seminorms.  Let $L \in
  \real^{n \times n}$ be the Laplacian matrix corresponding to an adjacency
  matrix $A \in \real^{n \times n}$ without self-loops, and let $d_{\rm
    out} = A \vect 1_n$. For each $i \in \{1, 2, \dots, n\}$, sort the
  off-diagonal entries of $A \vect e_j$ according to
  \[
  a_{(1),j} \ge a_{(2),j} \ge \cdots \ge a_{(n-1), j} .
  \]
  Then
  \begin{align*}
    \mu_{\dist,1}(-L) &= -\min_j \Bigg\{\![d_{\rm out}]_j\!-\!\!\!\sum_{i=1}^{\lfloor \frac n 2 \rfloor - 1} a_{(i), j}\!+\!\!\!\sum_{i = \lceil \frac n 2 \rceil}^{n-1} a_{(i), j}\!\Bigg\}, \\
    \mu_{\dist,2}(-L) &= \min_{b \in \real} \left\{
    b : \Pi_n L + L^\top \Pi_n \succeq -2 b \Pi_n
    \right\}, \\
    \mu_{\dist,\infty}(-L) &= -\min_{i \ne j} \Bigg\{\!a_{ij}\!+\!a_{ji}\!+\!\!\sum_{k\neq i,j} \min\{a_{ik}, a_{jk}\!\}
    \Bigg\}.
  \end{align*}  
\end{theorem}

\begin{proof}
  Set $S_h=I_n-hL$.  Observe that $S_h$ is row-stochastic for every $h>0$,
  and its entries are
  \begin{equation*}
    [S_h]_{ij} = \begin{cases}
      1 - h [d_{\rm out}]_i, & i = j, \\
      h a_{ij}, & i \ne j.
    \end{cases}
  \end{equation*}
  Also, $\mu_{\seminorm{\cdot}{}}(-L) = \lim_{h \to 0^+} h^{-1} \left(\seminorm{S_h}{} - 1 \right)$ for any seminorm $\seminorm{\cdot}{}$. 

  \noindent\emph{Case $\seminorm{\cdot}{\dist,1}$:} \quad
  For each $ j \in
  \{1, 2, \dots, n\}$, sort the entries of $S_h \vect e_j$ as
  \[
  (S_h)_{(1), j} \ge (S_h)_{(2), j} \ge \cdots \ge (S_h)_{(n), j}
  \]
  Assume $h$ is so small that $(S_h)_{(1), j} = (S_h)_{j, j}$. Then by~\eqref{eq:ind-norm-d1},
  \begin{multline*}
    \seminorm{S_h}{\dist,1} = \max_j \Bigg\{
    \sum_{i=1}^{\lfloor \frac n 2 \rfloor} s_{(i), j} - \sum_{i = \lceil \frac n 2 \rceil + 1}^n s_{(i), j}
    \Bigg\}  \\
    = 1 + h \max_j \Bigg\{ 
    - [d_{\rm out}]_j + \sum_{i=2}^{\lfloor \frac n 2 \rfloor} a_{(i), j} - \sum_{i = \lceil \frac n 2 \rceil + 1}^n a_{(i), j} \Bigg\}.
  \end{multline*}
  Substituting into \eqref{eq:log-norm-limit} yields the formula for
  $\mu_{\dist,1}(-L)$, since the order of the off-diagonal elements of $A
  \vect e_j$ is identical to the order of the off-diagonal elements of $S_h
  \vect e_j$ for all $h > 0$.

  \noindent\emph{Case $\seminorm{\cdot}{\dist,2}$:} \quad By
  \eqref{eq:ind-norm-d2},
  \begin{align*}
    &\seminorm{S_h}{\dist,2} = \min_{b \ge 0} \left\{
    b : S_h^\top \Pi_n S_h \preceq b^2 \Pi_n 
    \right\} \\
    &=  \min \left\{
    b\ge 0 : (I_n - h L)^\top \Pi_n (I_n - h L) \preceq b^2 \Pi_n 
    \right\}\\
    &=  \min  \left\{
    b \ge 0 : h^2 L^\top \Pi_n L - h \Pi_n L - h L^\top \Pi_n \preceq (b^2 - 1) \Pi_n 
    \right\}
  \end{align*}
  Therefore  $\frac{\seminorm{S_h}{\dist,2} - 1}{h}$ is equal to
  \begin{equation*}
    \min \{ \tfrac{b - 1}{h} : b \ge 0, h^2 L^\top \Pi_n L - h \Pi_n L - h
    L^\top \Pi_n \preceq (b^2 - 1) \Pi_n \}.
  \end{equation*}

  Let $\bar b = h^{-1}(b - 1)$, so that $b \ge 0$ if and only if $\bar b
  \ge -h^{-1}$, and $(b^2 - 1) = h \bar b (2 + h \bar b)$.  Performing this
  change of variables,
   \begin{multline*}
     \mu_{\dist,2}(L) = \lim_{h \to 0^+} \frac{\seminorm{S_h}{\dist,2} - 1}{h} \\
     =   \!\lim_{h \to 0^+}\!\! \min \{
     \bar b \ge -h^{-1} : h L^\top \Pi_n L - \Pi_n L - L^\top \Pi_n \preceq \bar b (2 + h \bar b) \Pi_n \} \\
        =  \min \! \left\{   \bar b \ge 0 : -\Pi_n L - L^\top \Pi_n \preceq 2 \bar b \Pi_n      \right\} ,
    \end{multline*}
   which is equivalent to the formula for $\mu_{\dist,2}(-L)$.

   \noindent\emph{Case $\seminorm{\cdot}{\dist,\infty}$:} \quad Assume $h$
   is sufficiently small that $1-h[d_{\rm out}]_i>ha_{ji}$ for all $i,
   j$. Applying \eqref{eq:ind-norm-dinf},
    \begin{align*}
      \seminorm{S_h}{\dist,\infty} &= 1 - \min_{i \ne j} \{ 
            \min\{1 - h [d_{\rm out}]_i, h a_{ji} \} \\
            & + \min\{1 - h [d_{\rm out}]_i, h a_{ij} \} + h\sum_{k \neq i, j} \min\{a_{ik}, a_{jk}\}
       \} \\
        &= 1 - h\min_{i \ne j} \left\{ 
            a_{ij} + a_{ji} + \sum_{k \neq i, j} \min\{a_{ik}, a_{jk}\}
        \right\} .
    \end{align*}
    Substituting into \eqref{eq:log-norm-limit} yields the formula for $\mu_{\dist,\infty}(-L)$.

\end{proof}

We also notice that, when $L=L^\top$, one can also show that
$\mu_{\dist,2}(-L)=-\lambda_2(L)$ (e.g., see \cite[Exercise~6.3]{FB:22}).
Explicit expressions for the $\ell_p$-projection logarithmic seminorm of
Laplacian matrices, for $p \in \{1,2,\infty\}$, are derived by
duality\begin{tac} and hence omitted\end{tac}.

\begin{arxiv}
\begin{corollary}[Explicit formulas for projection logarithmic seminorms] \label{thm:mu-explicit}
Consider the distance and the projection seminorms with ${\ker( \seminorm{\cdot}{\dist})}= {\ker( \seminorm{\cdot}{\proj})} = \spn\{\vect{1}_n\}$.    Let $L \in \real^{n \times n}$ be the Laplacian matrix corresponding to an adjacency matrix $A \in \real^{n \times n}$ without self-loops, and let $d_{\rm out} = A \vect 1_n$. For each $j \in \{1, 2, \dots, n\}$, sort the off-diagonal entries of $\vect e_i^\top A $ according to
$a_{i,(1)} \ge a_{i,(2)} \ge \cdots \ge a_{i,(n-1)}$
    Then
    \begin{align*}
          \mu_{\proj,1}(-L)\!\!&=\!\!-\min_{i \ne j} \left\{\!a_{ji}\!+\!a_{ij}\!+\!\!\sum_{k\neq i,j} \min\{a_{ki}, a_{kj}\!\}
        \right\}, \\
 \mu_{\proj,2}(-L) &= \min_{b \in \real} \left\{
            b : \Pi_n L + L^\top \Pi_n \succeq -2 b \Pi_n
        \right\}, \\
     \mu_{\proj,\infty}(-L)\!\!&=\!\!-\min_i \left\{\![d_{\rm out}]_i\!-\!\!\!\sum_{j=1}^{\lfloor \frac n 2 \rfloor - 1} a_{i, (j)}\!+\!\!\!\sum_{j = \lceil \frac n 2 \rceil}^{n-1} a_{i,(j)}\!\right\}.
    \end{align*}
\end{corollary}
\end{arxiv}

\section{Semicontracting dynamical systems}\label{sec:semicontraction}
We exploit now the duality result of induced matrix seminorms and induced matrix logarithmic seminorms for the study of strong semicontractivity of  dynamical systems. We also provide  some theoretical results that formalize semicontractivity conditions for linear and nonlinear dynamical systems both in discrete and continuous time.

    Given a vector subspace $\mcK \subset \real^n$ and a vector field $f: \real^n \to \real^n$, the \emph{perpendicular vector field} $f_\perp: \real^n \to \mcK^\perp$ and the \emph{parallel vector field} $f_\parallel: \real^n \to \mcK$ are denoted for all $x \in \real^n$  by $f_\perp(x) = \Pi_\perp f(x)$ and $f_\parallel(x) = (I_n - \Pi_\perp) f(x)$, respectively.  Given a seminorm $\map{\seminorm{\cdot}}{\real^n}{\realnonnegative}$, with kernel $\mcK$, the domain restriction of $\seminorm{\cdot}$ to $\mcK^\perp$, will be denoted  by $\map{\norm{\cdot}{\perp}}{\mathcal{K}^\perp}{\realnonnegative}$.

\begin{definition}[Invariant sets] \label{def:invariance}
    Let $f: \real^n \to \real^n$. A subspace $V \subset \real^n$ is \emph{$f$-invariant} on a domain $C \subseteq \real^n$ if $f(x + v) - f(x) \in V$ for all $x \in C$ and $v \in V$.
\end{definition}

\begin{lemma}[Differential characterization of invariance] \label{lem:invariance}
  Given a continuously differentiable map
  $\map{f}{C\subseteq\real^n}{\real^n}$, a subspace $V \subset \real^n$ is
  $f$-invariant if and only if $Df(x) V \subseteq V$ for all $x \in C$.
\end{lemma} 

\begin{proof}
    If $V$ is $f$-invariant, then $f(x + h v) - f(x) \in V$ for all $x \in C$, $v \in V$, and $h \in \real$, which implies that 
    \[
        Df(x)v = \lim_{h \to 0} \frac{f(x + h v) - f(x)}{h} \in V
    \]
    thus $Df(x) V \subseteq V$ for all $x \in C$. To prove the converse, assume $Df(x) V \subseteq V$; then for all $v \in V$,  
    \[
        f(x + v) - f(x) = \int_0^1 Df(x + \alpha v) v ~d\alpha \in V.
    \]
\end{proof}



\subsection{Discrete Time Semicontraction}

Let us consider the discrete time, time varying, nonlinear dynamics 
\begin{equation}\label{eq:dtsys}
    x(k+1) = f(k,x(k))
\end{equation}
with $k \in \integernonnegative, x \in \real^n$. We assume $f$ to be
continuously differentiable in the second argument.  In the following we
give a generalized definition of strongly semicontracting discrete time
system with respect to the one in \cite{SJ-PCV-FB:19q}. The generalization
applies to systems with arbitrary contraction step.

\begin{definition}[Semicontracting discrete time systems] \label{def:semicontracting_dt}
Let $\seminorm{\cdot}$ be a seminorm on $\real^n$ with kernel $\mathcal{K}.$ If there exists
$m \in \natural$, 
 $\rho<1$ and a domain $C \subseteq \real^n$ for which the time-varying vector field $f: \integernonnegative \times \real^n \rightarrow \real^n$ is such that
\footnote{$f^m$ is the $m$-the iterate of $f$ defined recursively by $f^m(k,x) = f(f^{m-1}(k,x))$.}
\begin{equation}\label{eq:weakened_def}
  \seminorm{D(f^m(k,x))} \leq \rho
\end{equation}
for all $ k \in \integernonnegative$ and  $x \in C$,
then the vector field is  \emph{strongly  semicontracting} on $C$ with rate
$\sqrt[m]{\rho}$.
\end{definition}



{\color{black} Lemma \ref{general_lem} provides sufficient conditions for
  two fundamental discrete-time systems to be strongly semicontracting.
\begin{lemma}[Strong semicontractivity of discrete-time affine systems]\label{general_lem}
Given a subspace $\mcK \subset \real^n$ and $p,q \in [1,\infty]$ with $p^{-1}+q^{-1} =1$, consider a sequence of matrices
$\{A(k)\}_{k\in\integernonnegative}\subset\real^{n\times{n}}$
satisfying:
 \begin{align}
   & A(k)\mcK \subseteq \mcK  \quad\text{for all } k \in \integernonnegative,
   \tag{invariance} \label{eq:invariance}
  \\
  &\rho\triangleq \sup_{k\in\integernonnegative} \tau_p(\mcK,A(k))<1. \tag{semicontractivity}
 \end{align}
\begin{enumerate}
    \item \label{invariance} Then  the system
\begin{equation}
    \label{sys:averaging_dt}
    x(k+1) = A(k)x(k)+b(k), \quad b(k)\in\real^n, 
\end{equation}
is strongly semicontracting with rate $\rho$ in the distance $\ell_q$ seminorm with kernel $\mcK$. Moreover
\begin{equation*}
\seminorm{x(k)-y(k)}{\dist,q} 
\leq \rho^k \seminorm{x(0)-y(0)}{\dist,q}.
\end{equation*}
\item \label{conservation}  The system
\begin{equation}
\label{sys:flow-system}
    x(k+1) = A^\top(k) x(k)+b(k), \quad b(k)\in\real^n,
\end{equation}
is strongly semicontracting with rate $\rho$ in the projection $\ell_p$ seminorm with kernel $\mcK$. Moreover, for any $x(0)$, $y(0)$ satisfying $x(0)-y(0)\in\mcK^\perp$,
\begin{equation*}
\seminorm {x(k)-y(k)}{\proj,p} \leq \rho^k \seminorm{x(0)-y(0)}{\proj,p}.
\end{equation*}
\end{enumerate}
\end{lemma}
\begin{proof}
The proof of part \ref{invariance} follows from equation \eqref{eq:ind-equiv-stoch} in Theorem \ref{mainth}, and from the conditional submultiplicative property $iii)$ Lemma \ref{lem:propertiesIMS}:
\begin{multline}\label{bound_averaging}
    \seminorm{x(k+1)-y(k+1)}{\dist,q}\\ \leq
  \seminorm{A(k)}{\dist,q} \seminorm{x(k)-y(k)}{\dist,q} \\
  = \tau_p(\mcK,A(k)) \seminorm{x(k)-y(k)}{\dist,q}.
\end{multline}
The proof of part \ref{conservation} follows from Lemma
\ref{lem:propertiesIMS} part \ref{prop1} since $x(k)-y(k) \in \mcK^\perp$,
$\forall k \in \integernonnegative$, as a consequence of the invariance
assumption~\eqref{eq:invariance} and therefore
\begin{multline}\label{bound_mc}
     \seminorm {x(k+1)-y(k+1)}{\proj,p}\\
     \leq \seminorm{A^\top(k)}{\proj,p} \seminorm{x(k)-y(k)}{\proj,p} \\
     = \tau_p(\mcK,A(k))\seminorm{x(k)-y(k)}{\proj,p}.
\end{multline}
\end{proof}
}

For example, when the subspace $\mcK$ is the consensus subspace, the matrices $\{A(k)\}_{k=0}^\infty$ are row-stochastic and the term $b(k) \equiv\vect{0}_n$ $\forall k \in \integernonnegative$, the systems~\eqref{sys:averaging_dt}
and~\eqref{sys:flow-system} are the standard averaging~\eqref{eq:A+Atopdt-A} and flow systems~\eqref{eq:A+Atopdt-AT} in the Introduction and the bounds \eqref{bound_averaging} and \eqref{bound_mc} are precisely the bounds~\eqref{eq:Markov-bound} and~\eqref{eq:averaging-bound}
stated in the Introduction.

The following theorem focuses on strong semicontractivity of discrete-time dynamical systems that enjoy the invariance property of the kernel of the seminorm.

\begin{theorem}[Discrete time semicontracting dynamics with invariance property]\label{mythe:contrperpdy}
Consider a system as in \eqref{eq:dtsys}. Let $\mcK \subset \real^n$ be an $f$-invariant subspace, and suppose that $f$ is strongly semicontracting with rate $\rho<1$, with respect to a seminorm $\seminorm{\cdot}{}$ on $\real^n$ with kernel $\mathcal K$. Then,
\begin{enumerate}
    \item \label{cascade_dt1} the system admits the cascade decomposition
\begin{align}
    x_\parallel (k+1) &=  f_{\parallel} (k,x_\parallel(k)+x_{\perp}(k)), \label{eq:dynamics-par} \\
    x_{\perp}(k+1) &= f_{\perp}(k,x_\perp (k)) \label{eq:dynamics-perp};
\end{align} 
\item\label{perpdyn_dt1} the perpendicular dynamics \eqref{eq:dynamics-perp} are strongly contracting on $\mcK^\perp$ with rate $\rho$, with respect to $\map{\norm{\cdot}{\perp}}{\mcK^\perp}{\realnonnegative}$; and
\item \label{contraction_dt1}for any two trajectories $x(k), y(k)$ of \eqref{eq:dtsys},
    \begin{equation*}
        \seminorm{x(k) - y(k)}{} \le \rho^k \seminorm{x(0) - y(0)}{}
    \end{equation*}
    for all $k \in \integernonnegative$.
\end{enumerate}
\end{theorem}
\begin{proof}
Regarding part \ref{cascade_dt1}, the  cascade decomposition \eqref{eq:dynamics-par}-\eqref{eq:dynamics-perp} follows from the observation that
\begin{align*}
    x_\perp(k+1) &= \Pi_\perp f(k, x_\parallel(k) + x_\perp(k)) \\
    &= \Pi_\perp f(k, x_\perp(k)) = f_\perp(k, x_\perp(k))
\end{align*}
where the second equality is due to the $f$-invariance of $\mcK$. 
Part \ref{perpdyn_dt1} follows from 
\begin{multline*}
     \max_{\substack{y \perp \mathcal{K}\\k\geq 0}} \seminorm{Df_\perp(k,y)} = \max_{\substack{y \perp \mathcal{K}\\k\geq 0}} \seminorm{\Pi_\perp Df(k,y)} \notag\\
     =\max_{\substack{y \perp \mathcal{K}\\k\geq 0}} \seminorm{Df(k,y)} \leq \max_{\substack{x \in \real^n\\k\geq 0}} \seminorm{Df(k,x)}\leq \rho
\end{multline*}
where the second equality follows from the fact that for a generic matrix $A$, $\seminorm{A} = \seminorm{\Pi_\perp A}$.
Part \ref{contraction_dt1} is a direct consequence of \ref{perpdyn_dt1}.
\end{proof}

The following theorem focuses on strong semicontractivity of discrete-time dynamical systems that enjoy the invariance property of the orthogonal complement of the kernel of the seminorm.

\begin{theorem}[Discrete time semicontracting dynamics with conservation property]\label{the:dt_perpcontraction}
Consider a system as in \eqref{eq:dtsys}. Let $\mathcal K \subset \real^n$ such that $\mathcal K^\perp$ is an $f$-invariant subspace. 
Let $\map{f}{\integernonnegative\times \real^n}{\real^n}$ be strongly semicontracting with rate $\rho<1$ with respect to a seminorm $\seminorm{\cdot}$ on $\real^n$ with kernel $\mcK$. Then,
\begin{enumerate}
\item \label{cascade_conserv_dt} the system admits the cascade decomposition
\begin{align}
    x_\parallel (k+1) &= f_{\parallel}(k,x_\parallel (k)), \label{dt_parallel}\\
    x_{\perp}(k+1) &= f_{\perp} (k,x_\parallel(k)+x_{\perp}(k))\label{dtperpdyn};
\end{align}
\item\label{th:dtsdecoupled1} for each $x_\parallel \in \mcK$, the vector field $x_\perp \mapsto f_\perp (k,x_\parallel + x_\perp)$ is strongly contracting with rate $\rho$, with respect to  $\map{\norm{\cdot}{\perp}}{\mcK^\perp}{\realnonnegative}$;
\item \label{th:dtsdecoupled3} if the map $x_\parallel \mapsto f_\perp(k, x_\parallel + x_\perp)$ is Lipschitz\footnote{That is, for all $x_\parallel, y_\parallel \in \mcK$, $z_\perp \in \mcK^\perp$, and $k \in \integernonnegative$, we have
    $||f_\perp(k, x_\parallel\!+\!z_\perp)\!-\!f_\perp(k, y_\parallel\!+\!z_\perp)||_\perp \le \ell d_{\mcK}(x_\parallel, y_\parallel)$.} with constant $\ell \in \real$ with respect to some metric $d_{\mcK}$ on $\mcK$, 
    then for any two trajectories $x(k), y(k)$ of \eqref{eq:dtsys}, satisfying $x(0)-y(0) \in \mcK^\perp$
    \begin{align*}\label{eq:iss_dt}
        \begin{split}
            &\seminorm{x(k+1) - y(k+1)}{} \\
            &\qquad\le \rho \seminorm{x(k) - y(k)}{} + \ell d_{\mcK}(x_\parallel(k), y_\parallel(k))
        \end{split}
    \end{align*}
    for all $k \in \integernonnegative$.
\end{enumerate}
 \end{theorem}
\begin{proof}

Regarding part \ref{cascade_conserv_dt}, the cascade decomposition \eqref{dt_parallel}--\eqref{dtperpdyn} follows from the observation that
\begin{align*}
    x_\parallel(k+1) &= (I_n - \Pi_\perp) f(k, x_\parallel(k) + x_\perp(k))\\
    &= (I_n - \Pi_\perp) f(k, x_\parallel(k)) = f_\parallel(k, x_\parallel(k))
\end{align*}
where the second equality is due to the $f$-invariance of $\mcK^\perp$. 
To prove \ref{th:dtsdecoupled1}, fix $x_\parallel \in \mcK$, and pick any $x_\perp, y_\perp \in \mcK^\perp$. Then
\begin{align*}
    &\norm{f_\perp(k, x_\parallel + x_\perp(k)) - f_\perp(k, x_\parallel + y_\perp(k))}{\perp} \\
    &\qquad= \seminorm{f_\perp(k, x_\parallel + x_\perp(k)) - f_\perp(k, x_\parallel + y_\perp(k))}{} \\
    &\qquad\le \rho \seminorm{x_\perp(k) - y_\perp(k)}{} = \rho \norm{x_\perp(k) - y_\perp(k)}{\perp}
\end{align*}
To prove \ref{th:dtsdecoupled3}, let $x_\parallel, y_\parallel \in \mcK$ and $x_\perp, y_\perp \in \mcK^\perp$. Then
\begin{align*}
    &\seminorm{f(k, x_\parallel(k) + x_\perp(k)) - f(k, y_\parallel(k) + y_\perp(k))} \\
    &\qquad \le \seminorm{f(k, x_\parallel(k) + x_\perp(k)) - f(k, y_\parallel(k) + x_\perp(k))} \\
    &\qquad + \seminorm{f(k, y_\parallel(k) + x_\perp(k)) - f(k, y_\parallel(k) + y_\perp(k))} \\
    &\qquad \le \ell d_\mcK(x_\parallel(k), y_\parallel(k)) + \rho \seminorm{x(k)- y(k)}
\end{align*}
where the first inequality is due to the subadditivity property and the second one follows from point \ref{th:dtsdecoupled1} and the invariance of $\mcK^\perp$.
\end{proof} 

\subsection{Continuous Time Semicontraction} 
Let us consider the continuous time, time varying, nonlinear dynamics 
\begin{equation}\label{sys:ct}
    \dot{x}(t) = f(t,x(t))
\end{equation}
with $t \in \realnonnegative, x \in \real^n$. We assume $f$ to be continuously differentiable in the second argument.
\begin{definition}[Semicontracting continuous time systems] Let $\seminorm{\cdot}$ be a seminorm on $\real^n$ with kernel $\mathcal{K}.$ The time-varying vector field $f: \realnonnegative \times \real^n \rightarrow \real^n$ is \emph{strongly  infinitesimally semicontracting} with rate $c>0$  on a domain $C \subseteq \real^n$ if $\forall t \in \realnonnegative$ and  $x \in C$,
\begin{equation*}
    \mu_{\seminorm{\cdot}}({Df(t,x)}) \leq -c.
\end{equation*}
\end{definition}
{\color{black}
Lemma \ref{general_lem_ct} provides sufficient conditions for two fundamental continuous time dynamical systems to be strongly infinitesimally semicontracting. 
\begin{lemma}[Strong semicontractivity of continuous-time affine systems]\label{general_lem_ct}
Given a subspace $\mcK \subset \real^n$ and $p,q \in [1,\infty]$ with
$p^{-1}+q^{-1} =1$, consider a sequence of matrices
$\{A(t)\}_{t\in\realnonnegative}\subset\real^{n\times{n}}$ satisfying:
 \begin{align*}
  & A(t)\mcK \subseteq \mcK  \quad\text{for all } t \in \realnonnegative,    \tag{invariance} \\
  & c\triangleq -\sup_{t\in\realnonnegative} \mu_{{\dist,p}}(A(t))>0.    \tag{semicontractivity}
 \end{align*}
\begin{enumerate}
    \item \label{invariance_ct} The system
\begin{equation*}
    \label{sys:averaging}
    \dot{x}(t) = A(t)x(t)+b(t), \quad b(t)\in\real^n,
\end{equation*}
is strongly infinitesimally semicontracting with rate $c$ in the distance $\ell_p$ seminorm with kernel $\mcK$,   moreover 

\begin{equation*}
\seminorm{x(t)-y(t)}{\dist,p} 
\leq e^{-ct}\seminorm{x(0)-y(0)}{\dist,p}, \, \forall t
\end{equation*}
\item \label{conservation_ct}   the system
\begin{equation}
\label{sys:flow-system-ct}
    \dot{x}(t) = A^\top(t) x(t)+b(t), \quad b(t)\in\real^n,
\end{equation}
is strongly infinitesimally semicontracting with rate $c$ in the projection $\ell_q$ seminorm with kernel $\mcK$, moreover, for any $x(0), y(0)$ satisfying $x(0)-y(0) \in \mcK^\perp$,
\begin{equation}\label{bound_flow_ct}
\seminorm {x(t)-y(t)}{\proj,q} \leq e^{-ct} \seminorm{x(0)-y(0)}{\proj,q}, \, \forall t.
\end{equation}

\end{enumerate}
\end{lemma}
\begin{proof}
The proof of part \ref{invariance_ct}
follows from Theorem 13, part i) in \cite{SJ-PCV-FB:19q}.
To prove part \ref{conservation_ct} we follow a similar reasoning as in Theorem 11  from \cite{SJ-PCV-FB:19q}. In fact, for all $x(0),y(0)$ such that $x(0)-y(0) \in \mcK^\perp$,  since the solutions $t \mapsto x(t)$ of \eqref{sys:flow-system-ct} are differentiable, by defining $z(t) \triangleq x(t)-y(t)$, for small $h$, one can write
\begin{align*}
   z(t+h) &= z(t)+
    h(A^\top(t) (z(t))) + o(h)= \Pi_\perp (z(t+h))
\end{align*}
since $z(t) \in \mcK^\perp$ and $A^\top \mcK^\perp \subseteq \mcK^\perp$ by hypothesis. Therefore, by Lemma \ref{lem:projection} and Lemma \ref{lem:propertiesIMS} part \ref{prop1}
\begin{align*}
    &\frac{\seminorm{z(t+h)}-\seminorm{z(t)}}{h} \leq\notag\\
    &\frac{\seminorm{I_n+hA^\top (t)}-1}{h} \seminorm{z(t)}+ \frac{o(h)}{h}. 
\end{align*}
Taking the limit as $h \rightarrow 0^+$, one gets $\frac{d}{dt}
\seminorm{z(t)} \leq \mu_{\seminorm{\cdot}}(A^\top(t)) \seminorm{z(t)}$.
Finally, from the Gr\"onwall comparison inequality (e.g.,
see~\cite[Exercise~2.1]{FB:22-CTDS})
\begin{multline*}
    \seminorm{x(t)-y(t)} \leq \\{\rm exp} \Big( \int_0^t \mu_{\seminorm{\cdot}}(A^\top(\tau)) d\tau \Big) \seminorm{x(0)-y(0)}.
\end{multline*}
Eq. \eqref{bound_flow_ct} follows from the fact that $\mu_{{\dist,p}}(A(t))=\mu_{{\proj,q}}(A^\top(t)) \leq -c$ for all $t$. 
\end{proof}
}

The following theorem focuses on strong infinitesimal semicontractivity of continuous-time dynamical systems that enjoy the invariance property of the kernel of the seminorm. This theorem extends Theorem 13 from \cite{SJ-PCV-FB:19q} through the formulation of a cascade decomposition and by establishing a strong contractivity property on the orthogonal complement to the seminorm kernel. 

\begin{theorem}[Continuous time semicontracting dynamics with invariance property, partially from \cite{SJ-PCV-FB:19q}]\label{mythe:ctcontrperpdy} Consider a system as in \eqref{sys:ct}. Let $\mcK \subset \real^n$ be an $f$-invariant subspace and suppose that $f$ is strongly infinitesimally semicontracting with rate $c>0$, with respect to a seminorm $\seminorm{\cdot}$ in $\real^n$ with kernel $\mcK$. Then,
\begin{enumerate}
    \item\label{cascade_ct1} the system admits the cascade decomposition 
\begin{align}
    \dot{x}_\parallel(t) &=  f_{\parallel} (t,x_\parallel(t)+x_{\perp}(t)),\\
    \dot{x}_{\perp}(t) &= f_{\perp}(t,x_\perp (t));\label{ct_perpdyn}
\end{align}
\item \label{perp_ct1} the perpendicular dynamics \eqref{ct_perpdyn} are strongly infinitesimally contracting on $\mcK^\perp$ with rate $c$, with respect to $\map{\norm{\cdot}{\perp}}{\mcK^\perp}{\realnonnegative}$;
\item \label{perp_ct_dyn1}for any two trajectories $x(t), y(t)$ of \eqref{sys:ct},
\begin{equation*}
    \seminorm{x(t)-y(t)} \leq e^{-ct} \seminorm{x(0)-y(0)}
\end{equation*}
for all $t \in \realnonnegative$. 
\end{enumerate}


\end{theorem}
\begin{proof}
Regarding part \ref{cascade_ct1}, the cascade decomposition is obtained by following the same reasoning as in Theorem \ref{mythe:contrperpdy}.
Part \ref{perp_ct1} follows from
\begin{align*}
     \mu_{\seminorm{\cdot}}(Df_\perp(t,y) )&=  \mu_{\seminorm{\cdot}}(\Pi_\perp Df(t,y))\notag\\
    & \leq \mu_{\seminorm{\cdot}}(Df(t,x))\leq -c
\end{align*}
where the first equality follows from the fact that for a generic matrix $A$, $\mu_{\seminorm{\cdot}}(A) = \mu_{\seminorm{\cdot}}(\Pi_\perp A)$. Part \ref{perp_ct_dyn1} is a direct consequence of part \ref{perp_ct1}.
\end{proof}
The following theorem focuses on strong semicontractivity of continuous-time dynamical systems that enjoy the invariance property of the orthogonal complement of the kernel of the seminorm.

\begin{theorem}[Continuous time semicontracting dynamics with conservation property]
Consider a system as in \eqref{sys:ct}. Let $\mcK \subset \real^n$  be such that $\mcK^\perp$ is an $f$-invariant subspace. Let $\map{f}{\realnonnegative\times \real^n}{\real^n}$ be strongly infinitesimally semicontracting with rate $c>0$ with respect to a seminorm  $\seminorm{\cdot}$ on $\real^n$ with kernel $\mcK$. Then,
\begin{enumerate}
    \item the system admits the cascade decomposition
    \begin{align}
    \dot{x}_\parallel (t) &= f_{\parallel}(t,x_\parallel (t)), \label{ct_parallel}\\
    \dot{x}_{\perp}(t) &= f_{\perp} (t,x_\parallel(t)+x_{\perp}(t));\label{perpdyn}
\end{align}
\item\label{th:dtsdecoupled1} for each $x_\parallel \in \mcK$, the vector field $x_\perp \mapsto f_\perp (t,x_\parallel + x_\perp)$ is strongly infinitesimally contracting with rate $c$, with respect to $\map{\norm{\cdot}{\perp}}{\mcK^\perp}{\realnonnegative}$;
\item \label{th:dtsdecoupled3} if the map $x_\parallel \mapsto f_\perp(t, x_\parallel + x_\perp)$ is Lipschitz\footnote{That is, for all $x_\parallel, y_\parallel \in \mcK$, $z_\perp \in \mcK^\perp$, and $ t \in \realnonnegative$, we have
    $||f_\perp(t, x_\parallel\!+\!z_\perp)\!-\!f_\perp(t, y_\parallel\!+\!z_\perp)||_\perp \le \ell d_{\mcK}(x_\parallel, y_\parallel)$.} continuous with constant $\ell \in \real$ with respect to some metric $d_{\mcK}$ on $\mcK$, 
 then for any two trajectories $x(t), y(t)$ of \eqref{sys:ct}, satisfying $x(0)-y(0) \in \mcK^\perp$
    \begin{align*}
       \begin{split}\label{eq:iss}
      & D^+\seminorm{x(t)-y(t)} \\
      &\qquad \leq -c \seminorm{x(t)-y(t)}+\ell d_{\mcK}(x_{\parallel}(t), y_{\parallel}(t))
    \end{split}
    \end{align*}
    for all $t \in \realnonnegative$, where $D^+(\cdot)$ indicates the upper right Dini derivative \cite[Section~2.1]{FB:22-CTDS}.
\end{enumerate}
\end{theorem}
\begin{proof}
The proof follows the same arguments as Theorem \ref{the:dt_perpcontraction} for discrete time systems.

\end{proof}

\section{Graph Theoretical Conditions for Semicontractivity}\label{sec:graphs}
We now provide graph theoretical conditions for the
systems~\eqref{eq:A+Atopdt} and \eqref{eq:A+Atopct} to be semicontracting
with respect to $\ell_p$ distance and projection seminorms, for $p
\in\{1,2,\infty\}$. For the discrete time case, the following conditions
are topological abstractions of algebraic conditions in
\cite{JL-SM-ASM-BDOA-CY:11, ICFI-TMS:11}. Lemma \ref{lem:averaging_ct} is
novel.

\begin{lemma}[Topological conditions for discrete-time averaging systems]\label{lem:averaging}
 The averaging system~\eqref{eq:A+Atopdt-A} $x(k+1)=Ax(k)$ with $A$
 row stochastic is strongly semicontracting in the
\begin{enumerate}
    \item\label{ave_graphd1} $\ell_1$ distance consensus seminorm if $A$ is
      doubly stochastic and $\mathcal{G}(A)$ is strongly connected and
      aperiodic;
    \item\label{ave_graphd2} $\ell_2$ distance consensus seminorm if $A$ is
      doubly stochastic and $\mathcal{G}(A)$ is weakly connected with self
      loops at each node;
    \item \label{ave_graphdinf} $\ell_\infty$ distance consensus seminorm
      if $\mathcal{G}(A)$ has self loops at each node and a globally
      reachable node.
\end{enumerate}
\end{lemma}
\begin{proof}
Condition \ref{ave_graphd1} ensures, in particular, that there exists $m \in \natural$ such that $A^m$ has at least $\lfloor \frac{n}{2}\rfloor+1$ nonzero entries in each column so the expression in \eqref{eq:ind-norm-d1} takes value less than one. Consequently, the system is strongly semicontracting according to condition \eqref{eq:weakened_def} in Definition \ref{def:semicontracting_dt}.

Condition \ref{ave_graphd2} directly follows from Lemma
\ref{prop:inequalities} and Theorem 8 in
\cite{JL-SM-ASM-BDOA-CY:11}. Finally, according to Corollary 4.5 in
\cite{FB:22}, condition \ref{ave_graphdinf} ensures that there exists $m
\in \natural$ such that $A^m$ (has a column with all nonzero entries and
hence) is scrambling. Consequently, according to Corollary 3.9 in
\cite{ICFI-TMS:11} and condition \eqref{eq:weakened_def} in Definition
\ref{def:semicontracting_dt} the system is strongly semicontracting.
\end{proof}

Strong semicontractivity of Markov chains in the $\ell_p$ projection
seminorms, $p \in \{1,2,\infty\}$, can be derived by duality.



\begin{lemma}[Topological conditions for continuous-time averaging]\label{lem:averaging_ct}
  The averaging system~\eqref{eq:A+Atopct} $\dot{x}=-Lx$ with $L$ the
  Laplacian of a graph with adjacency matrix $A$ and without self-loops, is
  strongly infinitesimally semicontracting in the
\begin{enumerate}
    \item\label{ct_dist1} $\ell_1$ distance consensus seminorm if $A$ is
      doubly stochastic and every node has at least $\lfloor
      \frac{n}{2}\rfloor$ in-neighbors,
   \item\label{ct_dist2} $\ell_2$ distance consensus seminorm if $A$ is
     doubly stochastic and $\mathcal{G}(A)$ is weakly connected,
    \item \label{ct_distinf} $\ell_\infty$ distance consensus seminorm if
      every two nodes are either (weakly) adjacent or have a common
      out-neighbor.
\end{enumerate}
\end{lemma}
\begin{proof}
To prove \ref{ct_dist1}, note that $\mu_{\seminorm{\cdot}{\dist,1}}(-L) <
0$ if and only if
\begin{equation*}
    \sum_{i=1}^{\lfloor \frac{n}{2} \rfloor-1} a_{(i),j}-\sum_{j = \lceil \frac{n}{2} \rceil}^{n-1} a_{(i),j} <1, \quad \forall i
\end{equation*}
that for $A$ doubly stochastic is fulfilled if and only if each node has at least $\lfloor \frac{n}{2}\rfloor$ in-neighbors.

To prove \ref{ct_dist2} note that for $A$ doubly stochastic $L\Pi_n = \Pi_n
L$ and hencethe formula for $\mu_{\dist,2}(-L)$ in
Theorem~\ref{thm:mu-explicit} reads as
\begin{equation*}
  \mu_{\dist,2}(-L)=\min_{b}\Big\{b: \frac{L+L^\top}{2}+b I_n \succeq 0 \text{ on } \mcK^\perp\Big\}.
\end{equation*}
The minimum is obtained for $b= -\lambda_2\Big(\frac{L+L^\top}{2} \Big)$
so that, for $\mathcal{G}(A)$ weakly connected, $\mu_{\dist,2}(-L)<0$.

To prove \ref{ct_distinf} note that
$\mu_{\seminorm{\cdot}{\dist,\infty}}(-L)<0$ if and only if
\begin{equation*}
  a_{ij}+a_{ji}+ \sum\nolimits_{k \neq i,j} \min\{a_{ik},a_{jk}\} >0 \quad \forall i \neq j
\end{equation*}
that for nonnegative adjacency matrices is true if and only if: $(i,j)$ is
an edge or, $(j,i)$ is an edge or $(i,k)$ and $(j,k)$ are an edge for some
third node $k$.
\end{proof}


\section{Conclusions}\label{sec:conclusions}
We have studied seminorms on vector spaces and induced matrix seminorms for
discrete- and continuous-time dynamical systems. We have shown how the
natural distance and projection seminorms are dual and how the long-studied
$\ell_p$ ergodic coefficients of a row-stochastic matrix are precisely
induced matrix seminorms. We have provided a comprehensive treatment of
semicontraction for discrete- and continuous-time systems with invariance
or conservation properties.
Future research directions include the application of semicontraction
theory to systems with symmetries, such as robotic vehicles ($SE(3)$
symmetry) and coupled oscillators (torus symmetry), as well as systems with
invariance properties, such as population games and evolutionary dynamics
(whose state space is the simplex).
A long-term elusive task is the definition of an ergodic coefficient that
is strictly less than unity for row-stochastic matrices satisfying weak
connectivity properties.


\appendices
\section{Seminorm Coefficients}\label{sec:appendix}

Here we recall some useful properties of standard $p$-norms. In the following, for a differentiable function $\map{f}{\real^n}{\real}$, we denote by $\nabla(f)$ its gradient. 

\begin{lemma}[Properties of differentiable $p$-norms \cite{WHY:90}] \label{lem:norms}
    Let $p \in (1, \infty)$, 
    then $||\cdot||_p$ has the following properties:
    \begin{enumerate}
        \item\label{norms:pt1} $||x||_p$ is differentiable on $\real^n$;
        \item\label{norms:pt2} $||x||_p = x^\top \nabla(\norm{x}{p})$ for all $x \in \real^n$; 
        \item \label{norms:pt3} $||\nabla(\norm{x}{p})||_q = 1$ for all $x \ne \vect 0_n$.
    \end{enumerate}
\end{lemma}   
\begin{proof}
See  the final remark and Equation (18) from \cite{WHY:90}.
\end{proof}

Based on Lemma \ref{lem:norms}, we  establish a novel and useful characterization of the distance and projection seminorms.

\begin{lemma}[Coefficients for distance seminorms] \label{lem:coefs}
    Let $p, q \in [1, \infty]$  be such that $p^{-1} + q^{-1} = 1$ and {\color{black} let $\mcK \subset \real^n$ be a vector subspace}. There exists a \emph{distance coefficient map} $\psi_p: \real^n \to \mcK^\bot$ such that, for all $x \in \real^n$,
    \begin{enumerate}
        \item $\psi_p(x) = \vect 0_n$ if $x \in \mcK$ and $\seminorm{\psi_p(x)}{\proj,q} = 1$ otherwise, and
        \item $\seminorm{x}{\dist,p} = \psi_p(x)^\top x$.
    \end{enumerate} 
\end{lemma} 

\begin{proof} 
    Let $V \in \real^{n \times k}$ be a a matrix whose columns are a basis for $\mcK$, so that we can write
    \[
        \seminorm{x}{\dist, p} = \min_{\alpha \in \real^k} \norm{x - V \alpha}{p}
    \]
    At the optimum $\alpha^*$, $\vect 0_k$ is a subgradient of $||x - V \alpha^*||_p$:
    \[
        \vect 0_n \in \partial ||x - V \alpha^*||_p = -V^\top G_p(x - V \alpha^*)
    \]
    where $G_p$ is the subdifferential $G_p = \partial ||\cdot||_p \subset \real^n$. Consequently, there exists a vector $\psi_p(x) \in G_p(x - V \alpha^*)$ such that $\psi_p(x) \in \ker(V^\top) = \mcK^\perp$. Note that $\seminorm{\psi_p(x)}{\proj, q} = \norm{\psi_p(x)}{q}$, that $\psi_p(x)^\top x = \psi_p(x)^\top (x - V \alpha^*)$, and that $\seminorm{x}{\dist, p} = \norm{x - V \alpha^*}{p}$, so we need only to show for each $p \in [1, \infty]$ that $\norm{\psi_p(x)}{ q} = 1$ and that $\psi_p(x)^\top (x - V \alpha^*) = \norm{x - V \alpha^*}{p}$.
    
    \paragraph*{Case $p = 1$}
    Using the standard formula for the subgradient of the absolute value function \cite{SB-LV:04}, $\psi_1(x) \in G_1(x - V \alpha^*)$ implies that
    \[
        (\psi_1(x))_i = \begin{cases}
            {\rm sgn}((x - V \alpha^*)_i) & (x - V \alpha^*)_i \neq 0 \\
            -1 \text{ or } +1
            , & (x - V \alpha^*)_i = 0
        \end{cases}, \;\; \forall i
    \]
    If $x \notin \mcK$, then $x - V \alpha^* \ne \vect 0_n$, so $\norm{\psi_1(x)}{\infty} = 1$. Furthermore,
    \begin{align*}
        (x - V \alpha^*)^\top \psi_1(x)
        &=\!\!\!\!\!\!\!\!\sum_{i:(x - V \alpha^*)_i \ne 0}\!\!\!\! (x - V \alpha^*)_i (\psi_1(x))_i \\
        &=\!\!\!\!\!\!\!\!\sum_{i:(x - V \alpha^*)_i \ne 0}\!\!\!\! (x - V \alpha^*)_i \sgn\left( (x - V \alpha^*)_i \right) \\
        &= \norm{x - V \alpha^*}.
    \end{align*} 
    \paragraph*{Case $p \in (1, \infty)$}
    If $p \in (1, \infty)$, then $\norm{\cdot}{p}$ is differentiable, so $G_p(z) = \nabla \norm{z}{p}$ for all $z \in \real^n$, and thus $\psi_p(x) = \nabla ||x - V \alpha^*||_p$ (where the gradient is taken with respect to $x - V \alpha^*$). If $x \notin \mcK^\perp$, then $x - V \alpha^* \ne \vect 0_n$, so $||\psi_p(x)||_q = 1$ due to Lemma \ref{lem:norms}. A further consequence of this lemma is that
    \[
        (x - V \alpha^*)^\top \psi_p(x) = ||x - V \alpha^*||_p = \seminorm{x}{\dist,p}
    \]
    
    \paragraph*{Case $p = \infty$}
    
    Let $\mathcal I \subseteq \{1, 2, \dots, n\}$ be the set of indices such that $\norm{x - V \alpha^*}{\infty} = |x - V \alpha^*|_i$. Using a standard formula for the subdifferential of a pointwise maximum \cite{SB-LV:04}, $\psi_\infty(x) \in G_\infty(x - V \alpha^*)$ implies that
    \[
        \psi_\infty(x) \in \conv \bigcup_{i \in \mathcal I} \partial |x - V \alpha^*|_i
    \]
    where $\conv$ denotes the convex hull, and the subdifferential of each absolute value is with respect to its argument. Therefore, there exist $g_i \in \partial |x - V \alpha^*|_i$ for each $i \in \mathcal I$, as well as convex weights $\lambda_i$, such that
    \[
        \psi_\infty(x) = \sum_{i \in \mathcal I} \lambda_i g_i
    \]
    For each $g_i$, we have $[g_i]_j = 0$ for $j \ne i$, since $|z|_i$ only depends on $z_i$ for any $z \in \real^n$. Furthermore, if $x \notin \mcK$, then $x - V \alpha^* \ne \vect 0_n$, so $|x - V \alpha^*|_i > 0$ for all $i \in \mathcal I$, which implies that $[g_i]_i = \sgn(x - V \alpha^*)_i$. Together, these two observations imply that 
    \[
        \norm{\psi_\infty(x)}{1} = \sum_{i \in \mathcal I} \lambda_i \norm{g_i}{1}
        = \sum_{i \in \mathcal I} \lambda_i = 1
    \]
    Finally,
    \begin{align*}
        &(x - V \alpha^*)^\top \psi_\infty(x)
        = \sum_{i \in \mathcal I} (x - V \alpha^*)_i \sum_{j \in \mathcal I} \lambda_j [g_j]_i \\
        &= \sum_{i \in \mathcal I} \lambda_i (x - V \alpha^*)_i \sgn(x - V \alpha^*)_i 
        = \norm{x - V \alpha^*}{\infty}
    \end{align*}
    
\end{proof}

\begin{lemma}[Coefficients for projection seminorms] \label{lem:coefs_disa}
    Let $p, q \in [1, \infty]$  be such that $p^{-1} + q^{-1} = 1$ and $\mcK \subset \real^n$ be a vector subspace. There exists a projection coefficient map $\zeta_p: \real^n \to \mcK^\perp$ such that, for all $x \in \real^n$,
    \begin{enumerate}
        \item $\zeta_p(x) = \vect 0_n$ if $x \in \mcK$ and $\seminorm{\zeta_p(x)}_{\dist,q} \leq 1$ otherwise, and
        \item $\seminorm{x}{\proj,p} = \zeta_p(x)^\top x$.
    \end{enumerate}
\end{lemma}

\begin{proof}  

Let $x \in \real^n$ and define $x_\perp = \Pi_\perp x$.

\paragraph*{Case $p = 1$} 
Let $\zeta_1(x) = \Pi_\perp\sgn(x_\perp)$. By Lemmas \ref{lem:projection} and \ref{lem:seminorm-props} \ref{prop:order-1},
\[
    \seminorm{\zeta_1(x)}{\dist, \infty} = \seminorm{\sgn(x_\perp)}{\dist, \infty} \le \norm{\sgn(x_\perp)}{\infty} \le 1,
\]
where $x \in \mcK$ implies that $\sgn(x_\perp) = \vect 0_n$. Furthermore,
\[
    \zeta_1(x)^\top x = \sgn(x_\perp)^\top \Pi_\perp x = \norm{x_\perp}{1} = \seminorm{x}{\proj, 1}.
\]

\paragraph*{Case $p \in (1, \infty)$}
Let $\zeta_p(x) = \Pi_\perp \nabla(\norm{x_\perp}{p})$. By Lemmas \ref{lem:projection}, \ref{lem:seminorm-props} \ref{prop:order-1}, and \ref{lem:norms} \ref{norms:pt3}, if $x \notin \mcK$, then
\[
    \seminorm{\zeta_p(x)}{\dist, q} = \seminorm{\nabla(\norm{x_\perp}{p})}{\dist, q} \le \norm{\nabla(\norm{x_\perp}{p})}{q} = 1.
\]
But if $x \in \mcK$, then $x_\perp = \vect 0_n$, so $\zeta_p(x) = \vect 0_n$. Furthermore, as a consequence of Lemma  \ref{lem:norms} \ref{norms:pt2}, 
\[
    \zeta_p(x)^\top x = (\nabla(\norm{x_\perp}{p}))^\top \Pi_\perp x = \norm{x_\perp}{p} = \seminorm{x}{\proj, p}.
\]

\paragraph*{Case $p = \infty$}

Let $i \in \{1, 2, \dots, n\}$ be such that $\norm{x_\perp}{\infty} = |x_\perp|_i$, and let $\zeta_\infty(x) = \sgn(x_\perp)_i\Pi_\perp \vect{e}_{i}$. By Lemmas \ref{lem:projection} and \ref{lem:seminorm-props} \ref{prop:order-1},
\[
    \seminorm{\zeta_\infty(x)}{\dist, 1}\!=\!\seminorm{\sgn(x_\perp)_i \vect e_i}{\dist, 1}
    \le \norm{\sgn(x_\perp)_i \vect e_i}{1} \le 1
\]
where $x \in \mcK$ implies that $\sgn(x_\perp) = \vect 0_n$. Furthermore,
\[
    \zeta_\infty(x)^\top x = \sgn(x_\perp)_i \vect e_i^\top x_\perp = \norm{x_\perp}{\infty} = \seminorm{x}{\proj, \infty}.
\] 

\end{proof}  

\section{Conjectures}
\begin{conjecture}[Optimal deflation and distance seminorm]\label{conj:distdefl}
Here a conjecture on the equivalence between $\ell_p$ distance seminorm and  $p$-optimal deflation as in Definition~\ref{def:optdef}.
    For each $p \in [1, \infty]$ and row-stochastic matrix $A \in \real^{n \times n}$,
    \begin{equation*}
        |A|_p = \seminorm{A}{\dist,p}. 
    \end{equation*}
\end{conjecture}
Here some reasons in support of this conjecture.
\begin{enumerate}
\item Expressions given in \cite{JL-ASM-BDOA-CY:11} for $p\in \{1,2,\infty\}$ of $x \in \real^n$ and of $A\in \real^{n\times n}$ row stochastic, of $|x|_p$ and $|A|_p$ coincide with the ones of $\seminorm{x}{\dist,p}$ and $\seminorm{A}{\dist,p}$, respectively.
\item If the envelope theorem \cite[Theorem~1.F.1]{AT:85} could be applied\footnote{It requires to prove the continuous differentiability of $v$ in \eqref{eq:optdef} as function of the inner optimization vector related to the induced matrix norm.} to the projection seminorm, with kernel $\mcK= \spn\{\vect{1}_n\}$, it would lead to the orthogonality constraint with respect to $\mcK$ and consequently to the equivalence between $|A|_p$ and $\seminorm{A^\top}{\proj,q}$.
\end{enumerate}

\begin{arxiv}
  \bibliographystyle{plainurl+isbn}
\end{arxiv}
\begin{tac}
  \bibliographystyle{IEEEtran}
\end{tac}
\bibliography{alias,Main,FB}

\begin{samepage}
\begin{IEEEbiography}[{\includegraphics[width=1in,height=1.25in,clip,keepaspectratio]{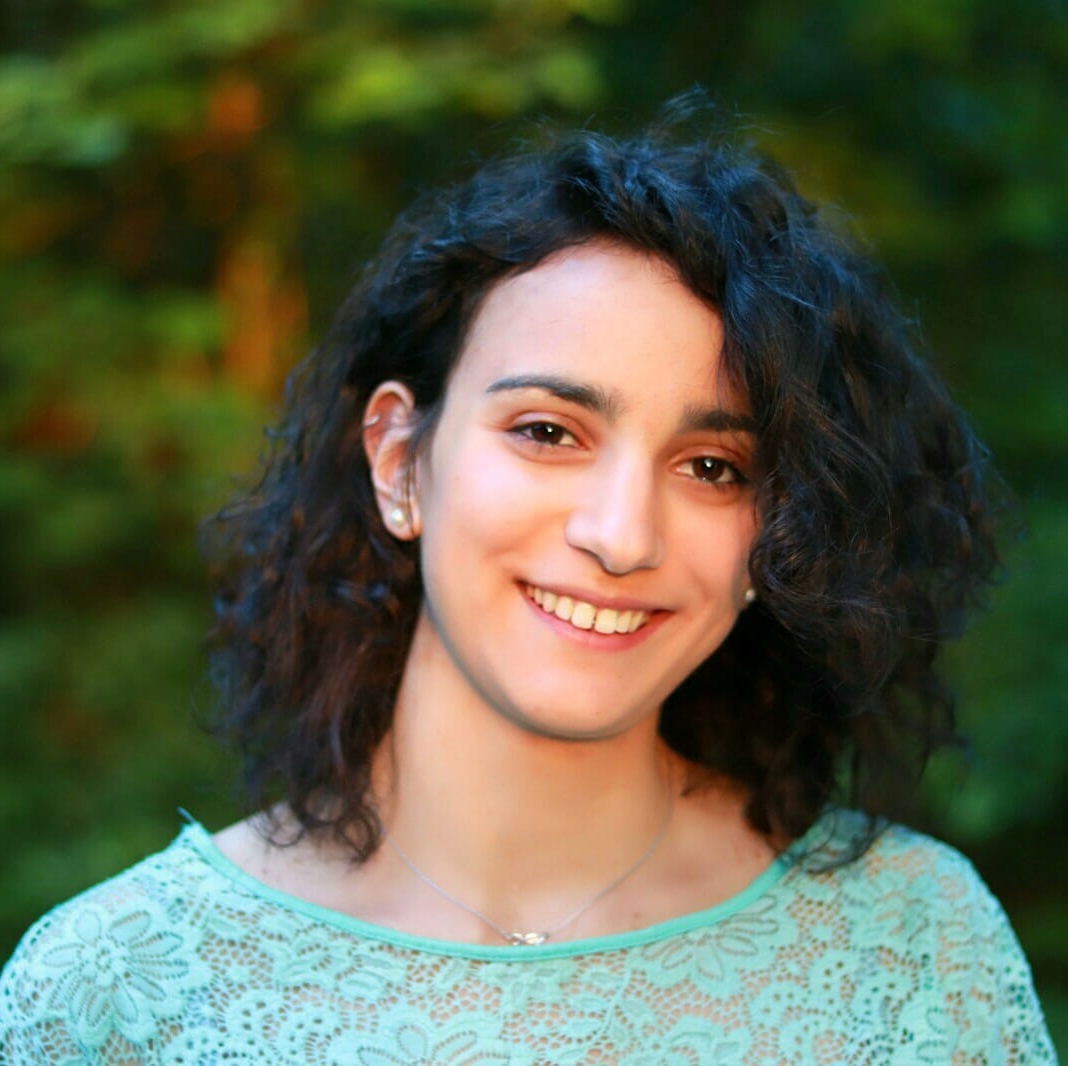}}]{Giulia De Pasquale} (S--'19) received  the  B.Sc. and the M.Sc  degree in Information Engineering and Systems and Control Engineering from the University of  Padova,  Italy,  in  2017  and  2019,  respectively, where  she  is  currently  pursuing  the  Ph.D.  degree with  the  Department  of  Information  Engineering. In 2022 she was visiting research scholar at the University of California, Santa Barbara. In 2018 and 2019 she was a visiting student at LTU, Sweden and ETH Z\"urich, respectively. Her current research interests are  in  the  analysis  and  control  of  networked systems, with a special focus on opinion dynamics.
\end{IEEEbiography}

\begin{IEEEbiography}[{\includegraphics[width=1in,height=1.25in,clip,keepaspectratio]{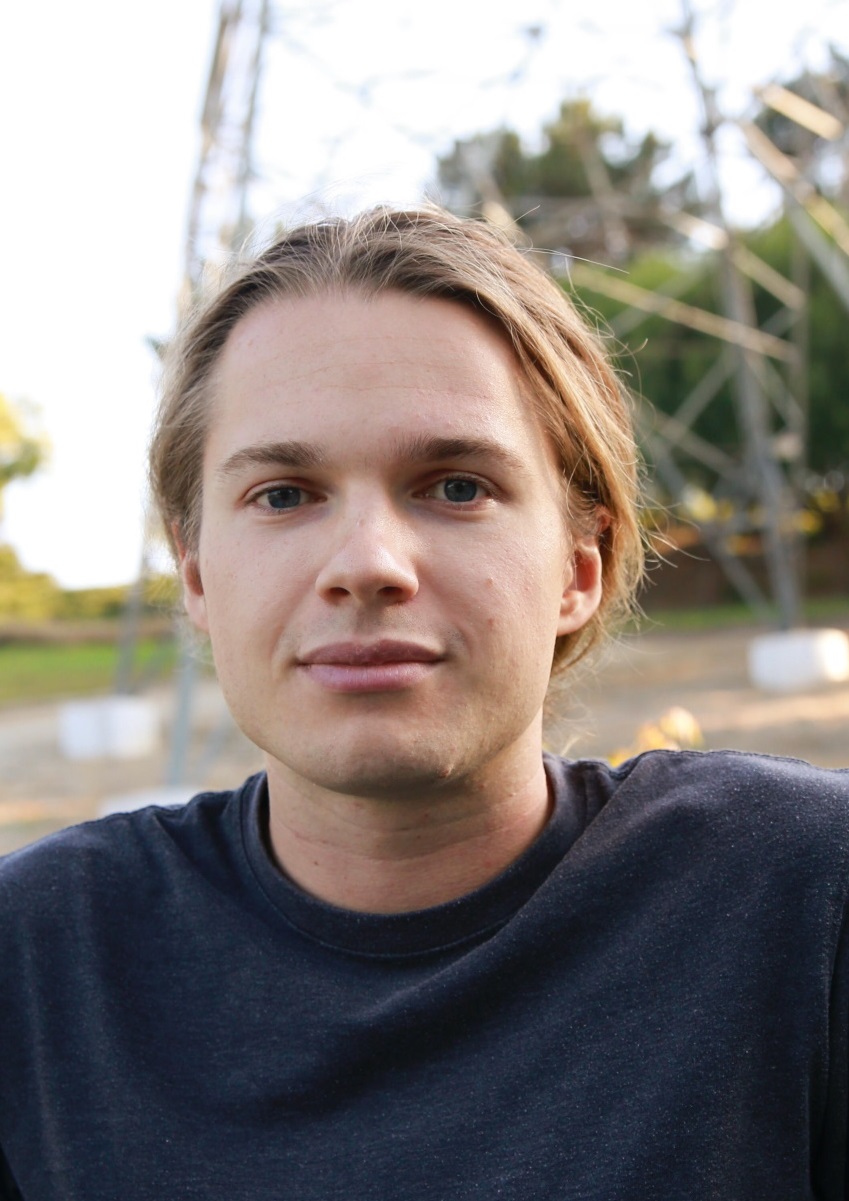}}]{Kevin D. Smith} (S--'16) received the B.S. degree in physics from Harvey Mudd College, Claremont, CA, USA, in 2017 and the M.S. degree in 2019 from the University of California, Santa Barbara, where he is currently working toward the Ph.D. degree with the Center for Control, Dynamical Systems, and Computation. His research interests include dynamics,
control, and identification of network systems, particularly infrastructure networks.
\end{IEEEbiography}

\begin{IEEEbiography}[{\includegraphics[width=1in,height=1.25in,clip,keepaspectratio]{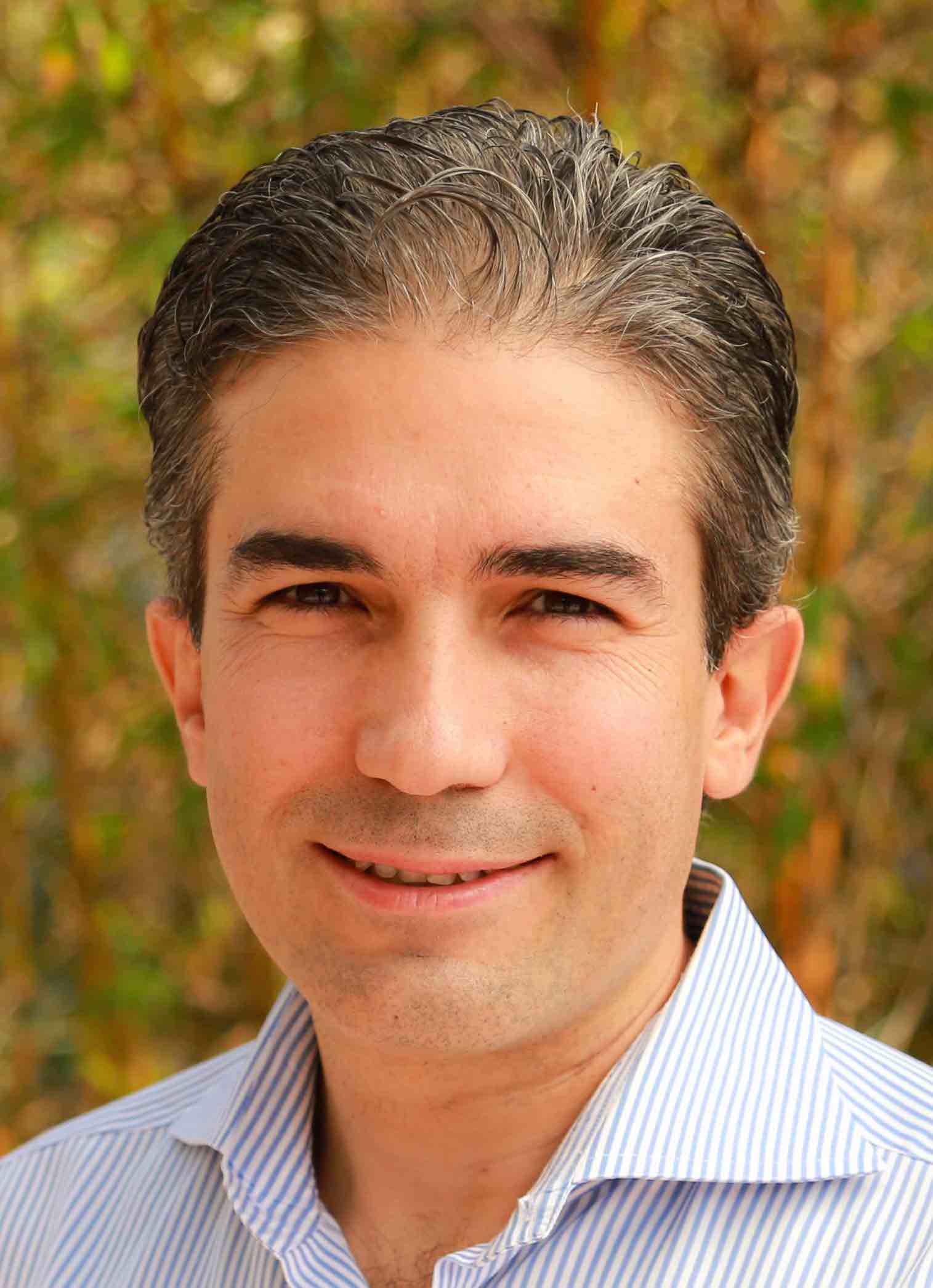}}]{Francesco Bullo} (S--’95 M--’99 SM--’03 F--'10)  is a Distinguished Professor of Mechanical Engineering with
the University of California, Santa Barbara, CA, USA. He was previously
with the University of Padova (Laurea degree, 1994), Italy, the California
Institute of Technology (Ph.D. degree, 1998), Pasadena, CA, and the
University of Illinois at Urbana-Champaign, IL, USA. His research interests
include contraction theory, network systems, and distributed control.  He
is the author or coauthor of Geometric Control of Mechanical Systems
(Springer, 2004), Distributed Control of Robotic Networks (Princeton,
2009), Lectures on Network Systems (KDP, 2022), and Contraction Theory for
Dynamical Systems (KDP, 2022, v1.0).  He served as IEEE CSS President and
SIAG CST Chair.  He is a Fellow of ASME, IEEE, IFAC, and SIAM.
\end{IEEEbiography}

\begin{IEEEbiography}[{\includegraphics[width=1in,height=1.25in,clip,keepaspectratio]{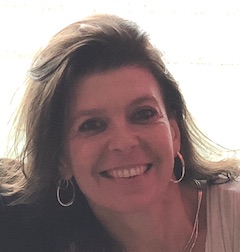}}]{M. Elena Valcher} (M--'99 SM--'03 F--'12) received the Master degree (1991) and the PhD degree (1995) from the University of Padova. Since  2005  she is full professor at the University of Padova. She is author/co-author of 87 journal papers, 105 conference papers, 2 text-books and several book chapters. Her research interests include social networks, cooperative control and consensus, positive switched systems and Boolean control networks. She is the Founding Editor in Chief of the IEEE Control Systems Letters (2017-). She was IEEE CSS President (2015). She received the 2011 IEEE CSS Distinguished Member Award and she is an IEEE Fellow since 2012.  She is a member of the IFAC Technical Board (2017-2020) and of the EUCA Board (2017-).
\end{IEEEbiography}

\end{samepage}

\end{document}